\pdfoutput=1
\documentclass[moor,nonblindrev]{informs1}              

\usepackage{natbib}
 \NatBibNumeric
 \def\bibfont{\small}%
 \bibpunct[, ]{[}{]}{,}{n}{}{,}%

\usepackage[colorlinks=true,breaklinks=true,bookmarks=true,urlcolor=blue,
     citecolor=blue,linkcolor=blue,bookmarksopen=false,draft=false]{hyperref}

\def\EMAIL#1{\href{mailto:#1}{#1}}
\def\URL#1{\href{#1}{#1}}         
\def\TheoremsNumberedThrough{\theoremstyle{TH}%
	\newtheorem{theorem}{\textbf{Theorem}} \newtheorem{lemma}{\textbf{Lemma}}
	     \newtheorem{assumption}{\textbf{Assumption}}
		\newtheorem{property}{\textbf{Property}}  \theoremstyle{EX}
	\newtheorem{remark}{\textbf{Remark}}   \newtheorem{definition}{\textbf{Definition}}   
}
\usepackage{algorithm}  
\usepackage{algorithmic}
\TheoremsNumberedThrough     

\EquationsNumberedThrough    

\begin{document}



\RUNTITLE{Global Convergence of IALM for 0/1-COP}

\TITLE{Global Convergence of Inexact Augmented Lagrangian Method for Zero-One Composite Optimization}

\ARTICLEAUTHORS{%
\AUTHOR{Penghe Zhang}
\AFF{Department of Applied Mathematics, Beijing Jiaotong University, Beijing, PR China, \EMAIL{19118011@bjtu.edu.cn} \URL{}}
\AUTHOR{Naihua Xiu}
\AFF{Department of Applied Mathematics, Beijing Jiaotong University, Beijing, PR China, \EMAIL{nhxiu@bjtu.edu.cn} \URL{}}
} 

\ABSTRACT{\textbf{Abstract.} We consider the problem of minimizing the sum of a smooth function and a composition of a zero-one loss function with a linear operator, namely zero-one composite optimization problem (0/1-COP). It is a versatile model including the support vector machine (SVM), multi-label classification (MLC), maximum rank correlation (MRC) and so on. However, due to the nonconvexity, discontinuity and NP-hardness of the 0/1-COP, it is intractable to design a globally convergent algorithm and the work attempting to solve it directly is scarce. In this paper, we first define and characterize the proximal stationarity to derive the minimum and the strongly exact penalization of the Lyapunov function, which is a variant of the augmented Lagrangian function for the 0/1-COP. Based on this, we propose an inexact augmented Lagrangian method (IALM) for solving 0/1-COP, where the subproblem is solved by the zero-one Bregman alternating linearized minimization (0/1-BALM) algorithm with low computational complexity. Under some suitable assumptions, we prove that the whole sequence generated by the IALM converges to the local minimizer of 0/1-COP. As a direct application, we obtain the global convergence of IALM under the assumption that the data matrix is full row rank for solving the SVM, MLC and MRC.}%



\maketitle

%
\section{Introduction.} In this paper, we consider the following zero-one composite optimization problem (0/1-COP):
\begin{equation} \label{eq1}
	\min_{w \in \mathbb{R}^p} f(w) + \lambda \| (Aw + b)_+ \|_0.
\end{equation}
By introducing an auxiliary variable $u$, problem (\ref{eq1}) can be equivalently reformulated as the following form
\begin{equation}  \label{eq2}
	\begin{aligned}
		&\min_{w \in \mathbb{R}^p, u \in \mathbb{R}^n} &&f(w) + \lambda \| u_+ \|_0,  \\
		&~~~~~~s.t. &&Aw + b = u,
	\end{aligned}
\end{equation}
where $f(w):\mathbb{R}^p \to \mathbb{R}$ is a smooth function, $\lambda$ is a positive weight parameter and $A \in \mathbb{R}^{n\times p}, b\in \mathbb{R}^n$. Given a vector $u = (u_1, \cdots, u_n)^\top \in \mathbb{R}^n$, we denote $u_+ := ( \max\{ 0,u_1 \},\cdots,\max\{ 0, u_n \} )^\top$. The $\| u \|_0$ is the $\ell_0$ norm of $u$, counting the number of its nonzero components. Therefore, $\| u_+ \|_0$ computes the number of all positive components of $u$. Moreover, $ \| u_+ \|_0 = \sum_{i = 1}^{n}\ell_{0/1} (u_i)$, where $\ell_{0/1}$ is called the zero-one loss function and defined as
\begin{equation} \notag
	\ell_{0/1} (u_i) : = \left\{  \begin{aligned}
		&    1,  &&u_i > 0, \\
		&    0,  &&u_i \leq 0.
	\end{aligned} \right.
\end{equation}   
The $\ell_{0/1}$ is an ideal function for recording the number of misclassification samples and indicating the rank correlation for a pair of samples. Therefore, (\ref{eq1}) is an original model for many problems in machine learning and statistics, including support vector machine (SVM), multi-label classification (MLC) and maximum rank correlation (MRC). Particularly, we give some examples in Section \ref{application} to exhibit the versatility of problem (\ref{eq1}). However, problem (\ref{eq1}) is nonconvex, discrete and NP-hard due to the $\| (\cdot)_+ \|_0$ term. Therefore, there are few studies on solving (\ref{eq1}) directly.



The augmented Lagrangian method (ALM) and its variants are well-known for solving constrained optimization problem. The ALM was first proposed independently by \citet{hestenes1969multiplier} and \citet{powell1969method} and was originally known as the method of multipliers. Since then, voluminous of theoretical and numerical developments have been witnessed over the past half-century, including the exact penalty property of augmented Lagrangian function (e.g. \citet{nocedal2006numerical,bertsekas2014constrained}) and the convergence properties of the ALM for solving convex problems (e.g. \citet{rockafellar1976augmented}). The algorithm has also been implemented in some efficient solvers for large-scale convex programming (see e.g. \citet{li2018qsdpnal,chen2016semismooth,yang2015sdpnal}). However, in the nonconvex and discontinuous setting, the exact penalty theory for augmented Lagrangian function is scarce and the global convergence properties of ALM is not well studied. 

\citet{chen2017augmented} considered a nonconvex non-Lipschitz constrained optimization problem and designed an inexact ALM with a non-monotone proximal gradient (NPG) method solving its subproblem. They showed that any accumulation point of the iteration sequence is a feasible point of original problem, namely the proposed algorithm has subsequence convergence property. Moreover, if the relaxed constant positive linear dependence (RCPLD) condition and the basic qualification (BQ) holds at such accumulation point, then it is a KKT point. \citet{song2020zero} proposed an inexact ALM with the similar framework to \citet{chen2017augmented} for solving the sum of a convex and twice continuously differentiable $f$, and a composite term $\| Aw \|_0$. Especially, they adopted an alternative scheme to solve the subproblem. They proved that each accumulation point of the iteration sequence is a KKT point under the sequence boundedness assumption. Furthermore, an accumulation point is a local minimizer of the original problem if $A$ is an orthogonal matrix. \citet{teng2020augmented} designed an augmented Lagrangian proximal alternating (ALPA) method for solving a sparse discrete constrained problem where the objective function is the sum of a smooth term, a nonsmooth term and a $\ell_0$ term. Specifically, the ALPA is an inexact ALM similar to that of \citet{chen2017augmented} with a proximal alternating minimization method (PALM) to solve the subproblem. Moreover, \citet{teng2020augmented} proved the subsequence convergence property of ALPA under the assumptions of Mangasarian-Fromovitz constraint qualification (MFCQ) and basic qualification (BQ). There are also some papers (e.g. \citet{cordova2021revisiting,jia2021augmented,kanzow2021augmented}) designing ALM with subsequence convergence property when solving optimization problems with disjunctive constraints, like cardinality or rank constraints. 

The alternating direction method of multipliers (ADMM) is an important variant of ALM and has been widely applied in machine learning, data analysis and signal processing. The studies of ADMM in the convex setting is comprehensive and clear, see, for instance, \citet{shefi2014rate,han2018linear,boyd}.  In recent years, there are some works addressing the convergence of ADMM in nonconvex setting. 
In \citet{li2015global}, a proximal ADMM (P-ADMM) is proposed for minimizing the sum of a smooth function and the composition of a lower semi-continuous (l.s.c.) function with a linear operator. The algorithm is convergent under the assumption that $A$ is full row rank and the objective function is semialgebraic and its smooth term has a bounded Hessian. For the same problem as that of \citet{li2015global}, \citet{boct2020proximal} establish the global convergence of a P-ADMM under the assumption that $A$ is full row rank and the objective function is semialgebric. \citet{bot2019proximal} considered a more general model than that of \citet{li2015global,boct2020proximal}, where the objective function has an additional l.s.c. term and the smooth term has two variables. In \citet{bot2019proximal}, a proximal linearized ADMM (PL-ADMM) is proposed and the global convergence is established under the assumption that $A$ is full row rank, the objective function is semialgebraic.  Some papers focus on solving the nonconvex composition problem with two variables and linear constraints which are general model including (\ref{eq1}). Most of these works established the convergence of ADMM-type algorithms under the assumption that the objective function satisfies Kurdyka-$\L$ojasiewicz (K$\L$) property. Furthermore, when it comes to the the global convergence analysis of the special case (\ref{eq1}), \citet{yashtini2020convergence,jia2021incremental} requires that $A$ is nonsingular. The convergence analysis for multi-block ADMM in nonconvex setting can refer to \citet{wang2019global,gao2020admm,wang2018con}. 

Recently, \citet{bolte2018nonconvex} proposed a novel generic algorithm framework named adaptive Lagrangian-based multiplier method (ALBUM), which includes an exact ALM and a P-ADMM, for minimizing the sum of a smooth function and a composition of a l.s.c. function with a general continuously differentiable operator $\mathcal{A}(w)$. The global convergence of ALBUM is constructed under the assumptions that the iteration sequence is bounded, $\mathcal{A}(w)$ is uniformly regular and the objective function is semialgebraic. Moreover, in the case that $\mathcal{A}(w)$ is linear, i.e. $\mathcal{A}(w) = A w$ for all $w \in \mathbb{R}^p$, \citet{bolte2018nonconvex} gave a PL-ADMM with full splitting form and low computational complexity. The similar algorithm also appears in \citet{boct2020proximal} and is derived by appropriate choices of proximal term. If the boundedness of the generated sequence and the semialgebraic property of objective function are satisfied, the global convergence for this kind of algorithm further requires that $A$ is nonsigular and the condition number of $A^\top A$ is less than 2.

It is noteworthy that \citet{bolte2018nonconvex} defined a Lyapunov function which is an augmented Lagrangain function with a proximal term, and gave the relationship of the critical points between the Lyapunov function and the original problem. The ALBUM proposed in \citet{bolte2018nonconvex} performs minimization on the Lyapunov function with respect to the primal variables, and the generated sequence ensures the sufficient decreasing of the Lyapunov function values, which facilitates the whole sequence convergence analysis. This brilliant technique also appears in \citet{wang2018con,li2015global,hong2016convergence}. The above facts lead us to consider the question: Is the Lyapunov function is an exact penalty function for 0/1-COP? To answer this question, we first define a proximal-type stationary point, which differs with the commonly used critical point or KKT point, to characterize the minimizer of the Lyapunov function with fixed Lagrange multiplier. The proximal-type stationary point is beneficial to algorithm design and the proximal behavior (see \citet[Theorem 1.25]{rockafellar1976augmented}) of proximal operator facilitates our convergence analysis. Then we show that a solution (proximal-type stationary point or minimizer) of 0/1-COP is also a solution of the Lyapunov function with fixed Lagrange multiplier if the penalty parameters are sufficiently large, namely the Lyapunov function is an exact penalty function for 0/1-COP. This result can be seen as an extension of the classical result for the smooth nonlinear programming (see e.g. \citet[Theorem 17.5]{nocedal2006numerical}). Moreover, we prove that a proximal-type stationary point of Lyapunov function with an appropriately selected Lagrange multiplier yields a proximal-type stationary point of 0/1-COP, which means that Lyapunov function is a strongly exact penalty function for 0/1-COP. 

Based on the (strongly) exact penalty property of the Lyapunov function and the technique for convergence analysis in \citet{bolte2018nonconvex}, we propose a novel inexact ALM (IALM) with global convergence for solving 0/1-COP. Under the assumption that $A$ is full row rank and $f$ is coercive, we show that the sequence generated by IALM is bounded and each of its accumulation point is a local minimizer of 0/1-COP. Additionally, if $f$ is strongly convex, the whole sequence converges to a local minimizer of 0/1-COP. The IALM algorithm generates a sequence of subproblems, the objective function of which is the sum of a smooth function with two primal variables and a l.s.c. $\| (\cdot)_+ \|_0$ term. This subproblem can be solved by the proximal alternating linearized minimization (PALM) method proposed by \citet{bolte2014proximal}. 
To further pursue a simple and explicit iteration formula, we utilize a Bregman distance function (see e.g. \citet{wang2018con}) and propose a 0/1-BALM algorithm with global convergence property and low computational complexity. Precisely speaking, 0/1-BALM is a variant of PALM focusing on the subproblems generated by IALM and has the similar computational complexity to the PL-ADMM in \citet{bolte2018nonconvex,boct2020proximal}, and 0/1-ADMM in \citet{wang2021support} when the parameters are properly selected, meanwhile its whole sequence convergence only requires that $A$ is full row rank. Overall, PL-ADMM and 0/1-ADMM could be seen as degenerate cases of our IALM equipped with 0/1-BALM which performs one step when solving the subproblems.

The main contributions of this paper are summarized as follows.
\begin{itemize}
	\item[$\bullet$] We define a proximal stationarity to derive the minimum and the (strongly) exact penalty properties of Lyapunov function for (\ref{eq1}). Precisely speaking, we first define a proximal stationary point to characterize the minimizer of Lyapunov function with a fixed Lagrange multiplier. Then for the Lyapunov exact penalty properties, we show that a solution (proximal-type stationary point or minimizer) of 0/1-COP is also a solution of the Lyapunov function with a fixed Lagrange multiplier if the penalty parameters are sufficiently large. For the Lyapunov strongly exact penalty properties, we prove that a proximal-type stationary point of Lyapunov function with an appropriately selected Lagrange multiplier yields a proximal-type stationary point of 0/1-COP.

	
	\item[$\bullet$] We propose an inexact ALM (IALM) with the global convergence property for solving (\ref{eq1}) under specific assumptions. If $f$ is coercive and $A$ is full row rank, then the sequence generated by IALM is bounded and each of its accumulation point is a local minimizer of 0/1-COP. If we further assume that $f$ is strongly convex, the whole sequence converges to a local minimizer of 0/1-COP.
	
	\item[$\bullet$] We design a zero-one Bregman alternating linearized minimization (0/1-BALM) algorithm for solving the subproblems generated by IALM with low computational complexity and weak assumption on $A$ ensuring global convergence. Moreover, our IALM equipped with 0/1-BALM performing one iteration degenerates into the PL-ADMM proposed by \citet{bolte2018nonconvex,boct2020proximal}, and 0/1-ADMM proposed by \citet{wang2021support} under specific parameter setting. 
	
	\item[$\bullet$] We present the specific forms of (\ref{eq1}) in the applications to SVM, MLC and MRC problems and obtain the global convergence of IALM under the assumption that $A$ is full row rank for these cases.   
\end{itemize}

The rest of this paper is organized as follows. We study the minimization of Lyapunov function with a fixed Lagrange multiplier in Section \ref{lya_min}. The (strongly) exact penalty properties of the Lyapunov function are presented in Section \ref{lya_ep}. We give the framework of IALM and 0/1-BALM in Section \ref{ialm}. Global convergence of the IALM is analyzed in Section \ref{convergence_analysis}. Some examples and algorithmic analysis of (\ref{eq1}) are presented in Section \ref{application}.

\section{Notations and Preliminaries.} For convenience, we define some notation employed throughout this paper. Given a subset $S \subseteq [n] := \{1,2,\cdots,n\}$, its complementary set and cardinality are $S^c$ and $|S|$ respectively. Given vector $u \in \mathbb{R}^n$ (resp. $A \in \mathbb{R}^{n \times p}$), its subvector (resp. submatrix) composed of entries (resp. rows) indexed on $S$ is denoted by $u_S$ (resp. $A_S$). The neighborhood of $w^* \in \mathbb{R}^p$ with radius $\delta > 0$ is denoted by $\mathcal{N}(w^*;\delta) := \{ w \in \mathbb{R}^p : \| w - w^* \| < \delta \}$. The combination of two vectors is denoted by $(w;u) = (w^\top u^\top)^\top$. Let $\textbf{1}_n \in \mathbb{R}^n$ be a vector with all entries being one. $\mathbb{N}$ ($\mathbb{N}^+$) denotes the set including all (positive) natural numbers. 
The symbol $\odot$ represents the Hadamard product. Given a symmetric matrix $Q \in \mathbb{R}^{n\times n}$, $Q \succ (\succeq) 0$ represents that $Q$ is positive definite (semidefinite). Given a set $\Theta \subseteq \mathbb{R}^n$, the distance function associated with $\Theta$ is defined by $ {\rm dist}(u;\Theta):=\inf_{\theta \in \Theta} \| u - \theta \| $.

Recall from \citet[Definition 8.3]{RockWets98} that for a proper, lower semi-continuous (l.s.c.) function $g: \mathbb{R}^n \to ( - \infty, \infty ]$ and a point $u$ with finite $g(u)$,  the limiting subdifferential is defined as
\begin{equation} \notag
	\partial g(u) := \{ v : \exists u^t \to u, v^t \to v ~{\rm with}~ \liminf_{z \to u^t} \frac{g(z) - g(u^t) - \langle v^t, z - u^t \rangle}{\| z - u^t \|} \geq 0 ~{\rm for}~{\rm each}~t \}.
\end{equation}
Recall from \citet[Definition 1.22]{RockWets98} that for a proper, l.s.c. function $g:\mathbb{R}^n \to (-\infty, \infty]$ and constant $\alpha > 0$, the proximal operator of function $g$ with constant $\alpha > 0$ is defined by
\begin{equation} \notag
	{\rm Prox}_{\alpha g} (z) : = \mathop{\arg\, \min}\limits_{u \in \mathbb{R}^n }~ g(u) + \frac{1}{2\alpha} \| u - z \|^2.
\end{equation}

Particularly, $\|\cdot\|_0$ is proper and l.s.c., and $(\cdot)_+$ is proper and continuous. Hence, the composition $\| (\cdot)_+ \|_0$ is proper and l.s.c. function. The computing formula for the limiting subdifferential and proximal operator of $\| (\cdot)_+ \|_0$ are given in \citet{wang2021support,zhou2021quadratic}.
\begin{property}[\citet{wang2021support,zhou2021quadratic}] \label{prop+0}
	Some properties of $\|(\cdot)_+\|_0$ are summarized as follows.
	
	(i) Given a vector $u \in \mathbb{R}^n$, the limiting subdifferential of $\|(\cdot)_+\|_0$ can be represented as
	\begin{equation} \label{sub+0}
		\partial \| u_+ \|_0 = \left\{ z \in \mathbb{R}^n : z_i \left\{ \begin{aligned}
			&=0,  &&{\rm if}~ u_i \neq 0, \\
			&\geq0, &&{\rm otherwise},
		\end{aligned} \right. i\in [n]\right\}.
	\end{equation}
	
	
	(ii) Given a vector $z \in \mathbb{R}^n$ and constant $\alpha >0$, the proximal operator of $\lambda \| (\cdot)_+ \|_0$ with parameter $\alpha$ can be represented as
	\begin{equation}\label{prox0+}
		\big[ 	{\rm Prox_{\alpha\lambda \|(\cdot)_+\|_0}} ( z )  \big]_i = {\rm Prox_{\alpha\lambda \|(\cdot)_+\|_0}} ( z_i ) = \left\{  \begin{aligned}
			& 0, && z_i \in (0,\sqrt{2\lambda\alpha}), \\
			&\{ 0,z_i \}, && z_i \in \{ 0, \sqrt{2\lambda \alpha} \}, \\
			&z_i, && z_i \in (-\infty,0)\cup (\sqrt{2\lambda\alpha},\infty),
		\end{aligned} \right. i \in [n].
	\end{equation}
\end{property}

Suppose that $g$ is a $l_g$-smooth function, namely for any $u, \overline{u} \in \mathbb{R}^n$, 
\begin{equation}\notag
	\| \nabla g(u) - \nabla g(\overline{u}) \| \leq l_g \| u - \overline{u} \|,
\end{equation}
then the following descent lemma holds (see e.g. \citet[Lemma 5.7]{beck2017first}).
\begin{lemma}[Descent Lemma (\citet{beck2017first})]\label{dl}Let $g: \mathbb{R}^n \to (- \infty, \infty]$ be a $l_g$-smooth function. Then for any $u$ and $\overline{u} \in \mathbb{R}^n$,
	\begin{equation}
		g(u) \leq g(\overline{u}) + \langle \nabla g(\overline{u}), u - \overline{u} \rangle + \frac{l_g}{2} \| u - \overline{u} \|^2. \notag
	\end{equation} 
\end{lemma}

The augmented Lagrangian and Lyapunov function of (\ref{eq2}) can be written as (\ref{eq3}) and (\ref{eq4}) respectively:
\begin{align} 
	&\mathcal{L}_\rho (w,u,z) := f(w) + \lambda \| u_+ \|_0 + \langle z, Aw + b - u \rangle + \frac{\rho}{2} \| Aw + b - u \|^2,   \label{eq3} \\
	&\mathcal{V}_{\rho,\mu} (w,u,z,v) := \mathcal{L}_\rho (w,u,z) + \frac{\mu}{2} \| w - v \|^2,  \label{eq4}
\end{align}  
where $\rho > 0$ and $\mu > 0$ are penalty parameters. Moreover, we define $(w;u)$, $z$ and $v$ as primal variable, Lagrange multiplier and Lyapunov variable respectively.

\section{Lyapunov Function Minimum.} \label{lya_min}
In this section, we will conduct optimality analysis on the following problem
\begin{equation} \label{exact_penalty}
	\min_{(w;u;v)\in \mathbb{R}^{p} \times \mathbb{R}^{n} \times \mathbb{R}^{p}} \mathcal{V}_{\rho,\mu} ( w,u,\widetilde{z},v ),
\end{equation}
where $\widetilde{z} \in \mathbb{R}^n$ is a fixed vector. Problem (\ref{exact_penalty}) is important for establishing the exact penalty theory and the algorithm design. Let us define a continuously differentiable function
\begin{equation} \notag
	\widetilde{h}(w,u,v):= f(w) + \frac{\mu}{2} \| w - v \|^2 +  \langle \widetilde{z}, Aw + b - u \rangle + \frac{\rho}{2} \| Aw + b - u \|^2.
\end{equation}
Then the objective function in (\ref{exact_penalty}) can be represented as 
\begin{equation} \notag
	\mathcal{V}_{\rho,\mu} ( w,u,\widetilde{z},v ) = \widetilde{h}(w,u,v) + \lambda \| u_+ \|_0.
\end{equation}
We first show some properties of function $\widetilde{h}$.
\begin{property}\label{prop_htilde}The function $\widetilde{h}$ has the following properties:
	
	(i) The gradient of $\widetilde{h}$ can be represented as
	\begin{equation} 
		\begin{aligned} \label{nabla_h3}
			\nabla \widetilde{h}(w,u,v) &= \left(  \begin{array}{c}
				\nabla_w \widetilde{h}(w,u,v) \\
				\nabla_u \widetilde{h}(w,u,v) \\
				\nabla_v \widetilde{h}(w,u,v)
			\end{array}  \right)  = \left(  \begin{array}{c}
				\nabla f(w) + \mu(w - v) + \rho A^\top \big( Aw + b - u  + \widetilde{z} / \rho \big) \\
				\rho(  u - Aw - b  - \widetilde{z} / \rho )\\
				\mu( v - w )
			\end{array}  \right).
		\end{aligned}
	\end{equation}
	
	(ii) If $f$ is convex, then $\widetilde{h}$ is convex. Moreover, if $f$ is $\sigma_f$-strongly convex, then $\widetilde{h}$ is $\sigma$-strongly convex, which means that for any $(w;u;v)$ and $(\overline{w};\overline{u};\overline{v})\in \mathbb{R}^{p} \times \mathbb{R}^{n} \times \mathbb{R}^{p}$, it holds that
	\begin{align} 
		\widetilde{h}(w,u,v) - \widetilde{h}(\overline{w},\overline{u},\overline{v}) \geq & \langle \nabla_w \widetilde{h}(\overline{w},\overline{u},\overline{v}), w - \overline{w}  \rangle + \langle \nabla_u \widetilde{h}(\overline{w},\overline{u},\overline{v}), u - \overline{u}  \rangle + \langle \nabla_v \widetilde{h}(\overline{w},\overline{u},\overline{v}), v - \overline{v}  \rangle \notag \\
		&+ \frac{\sigma}{2} (\| w - \overline{w} \|^2 + \| u - \overline{u} \|^2 + \| v - \overline{v} \|^2 ), \label{h_str_con}
	\end{align}
	where
	$\sigma:= \sigma_f \rho \mu/ ( \sigma_f(\mu + \rho ) + \rho\mu \| A \|^2 + 2\rho \mu )  $.
\end{property}

\proof{Proof.}
Since (i) is easy to be verify, we only prove (ii). First, suppose that $f$ is convex, let us prove that $\widetilde{h}$ is convex. We define 
\begin{equation}
	\widetilde{\phi}( w,u,v ): = \frac{\mu}{2} \| w - v \|^2 +  \langle \widetilde{z}, Aw + b - u \rangle + \frac{\rho}{2} \| Aw + b - u \|^2. \notag
\end{equation}
$\widetilde{\phi}$ is a quadratic function with Hessian 
\begin{equation} \label{hess_phi}
	\nabla^2 \widetilde{\phi}(w,u,v) = \left[ \begin{array}{ccc}
		\mu I + \rho A^\top A~   &-\rho A^\top  & -\mu I \\
		-\rho A  &\rho I & 0 \\
		-\mu I & 0 &\mu I
	\end{array} \right] = R_1 \Lambda_1 R_1^\top, 
\end{equation}
where 
\begin{equation}
	R_1 = \left[ \begin{array}{ccc}
		I   &-A^\top  & -I \\
		0  & I & 0 \\
		0 & 0 & I
	\end{array} \right],~\Lambda_1 = \left[ \begin{array}{ccc}
		0   &0  & 0 \\
		0  &\rho I & 0 \\
		0 & 0 &\mu I
	\end{array} \right]. \notag
\end{equation}
Since $R_1$ is nonsingular and $\Lambda_1$ is positive semidefinite, we know that $\nabla^2 \widetilde{\phi}(w,u,v)$ is positive semidefinite. It follows from $\widetilde{h} = f + \phi$ and the convexity of $f$ that $\widetilde{h}$ is convex.

Next, if $f$ is $\sigma_f$-strongly convex, let us prove that $\widetilde{h}$ is $\sigma$-strongly convex. It follows from (\ref{hess_phi}) that
\begin{align}
	&\nabla^2 (\widetilde{\phi}(w,u,v) + \frac{\sigma_f}{2} \|w\|^2 ) - \nabla^2 \frac{\sigma}{2} ( \| w \|^2 + \| u \|^2 + \| v \|^2 ) \notag \\
	= & \left[ \begin{array}{ccc}
		\mu I + \rho A^\top A + (\sigma_f - \sigma) I   &-\rho A^\top  & -\mu I \\
		-\rho A  &(\rho - \sigma) I & 0 \\
		-\mu I & 0 &(\mu - \sigma) I
	\end{array} \right] = R_2 \Lambda_2 R_2^\top, \notag
\end{align}
where
\begin{align}
	&	R_2 = \left[ \begin{array}{ccc}
		I   &-\frac{\rho}{\rho - \sigma}A^\top  & -\frac{\mu}{\mu - \sigma}I \\
		0  & I & 0 \\
		0 & 0 & I
	\end{array} \right], ~\Lambda_2 = \left[ \begin{array}{ccc}
		(\rho - \frac{\rho^2}{\rho - \sigma}) A^\top A + ( \sigma_f + \mu - \sigma - \frac{\mu^2}{\mu - \sigma} )I   &0  & 0 \\
		0  &(\rho - \sigma) I & 0 \\
		0 & 0 &( \mu - \sigma ) I
	\end{array} \right]. \notag
\end{align}
By the definition of $\sigma$ and direct computation, we can obtain
\begin{equation}
	(\rho - \frac{\rho^2}{\rho - \sigma}) \| A \|^2 + \sigma_f + \mu - \sigma - \frac{\mu^2}{\mu - \sigma}  > 0,~\rho - \sigma > 0,~\mu - \sigma > 0, \notag
\end{equation}
which means that $\Lambda_2$ is positive definite and hence $\nabla^2 \widetilde{\phi}(w,u,v)+ (\sigma_f/2) \|w\|^2 - \nabla^2 (\sigma/2) ( \| w \|^2 + \| u \|^2 + \| v \|^2 )$ is positive definite. This implies that $\widetilde{\phi}(w,u,v) + (\sigma_f/2) \|w\|^2$ is $\sigma$-strongly convex. Then for any $(w;u;v)$ and $(\overline{w};\overline{u};\overline{v})\in \mathbb{R}^{p} \times \mathbb{R}^{n} \times \mathbb{R}^{p}$, it holds that
\begin{align} 
	\widetilde{\phi}(w,u,v) - \widetilde{\phi}(\overline{w},\overline{u},\overline{v}) \geq & \langle \nabla_w \widetilde{\phi}(\overline{w},\overline{u},\overline{v}), w - \overline{w}  \rangle + \langle \nabla_u \widetilde{\phi}(\overline{w},\overline{u},\overline{v}), u - \overline{u}  \rangle + \langle \nabla_v \widetilde{\phi}(\overline{w},\overline{u},\overline{v}), v - \overline{v}  \rangle \notag \\
	& + \frac{\sigma - \sigma_f}{2} \| w - \overline{w} \|^2 + \frac{\sigma}{2} ( \| u - \overline{u} \|^2  + \| v - \overline{v} \|^2 ). \label{sc_phi}
\end{align}
Moreover, the $\sigma_f$-strong convexity of $f$ yields that
\begin{align}
	f (w) - f(\overline{w}) \geq \langle \nabla f(\overline{w}), w - \overline{w} \rangle + \frac{\sigma_f}{2} \| w - \overline{w} \|^2. \label{sc_f}
\end{align}
Adding (\ref{sc_phi}) and (\ref{sc_f}), we can obtain the $\sigma$-strong convexity of $\widetilde{h}$.
\hfill \Halmos
\endproof






Next lemma gives a sufficient condition ensuring the existence of solution for (\ref{exact_penalty}).

\begin{lemma}[Existence of the Solution]\label{existence}The problem (\ref{exact_penalty}) exists a global minimizer if  $f$ is coercive, which means that $\lim_{ \| w \| \to \infty } f(w) = \infty $.
\end{lemma}
\proof{Proof.} Let us first prove the coercivity of $\mathcal{V}_{\rho,\mu} (w,u,\widetilde{z},v)$. 
From the definition of $\widetilde{\phi}$, we have 
\begin{align} \notag
	\widetilde{\phi}( w,u,v ) = \frac{\mu}{2} \| w - v \|^2  + \frac{\rho}{2} \| Aw + b - u + \widetilde{z} / \rho \|^2 - \frac{1}{2\rho} \| \widetilde{z} \|^2 .
\end{align}
This implies that $\widetilde{\phi}$ is bounded below. Moreover, from the definition of Lyapunov function, we have 
\begin{align} \notag
	\mathcal{V}_{\rho,\mu} (w,u,\widetilde{z},v) = \widetilde{\phi}(w,u,v) + f(w) + \lambda \| u_+ \|_0.
\end{align}
Then $\mathcal{V}_{\rho,\mu} (w,u,\widetilde{z},v)$ is the sum of a coercive function and two bounded  below functions and hence it is coercive.

Fixing $(\overline{w};\overline{u};\overline{v}) \in \mathbb{R}^{p} \times \mathbb{R}^{n} \times \mathbb{R}^{p}$ and taking $M = \mathcal{V}_{\rho,\mu} (\overline{w},\overline{u},\widetilde{z},\overline{v})$, then the level set $\Omega_M := \{ (w;u;v): \mathcal{V}_{\rho,\mu} (w,u,\widetilde{z},v) \leq M \}$ is nonempty. Moreover, we claim that $\Omega_M$ is bounded. Indeed, if $\Omega_M$ is unbounded, then there exists a sequence $\{ (w^k;u^k;v^k) \}_{k \in \mathbb{N}} \subseteq \Omega_M$ satisfying $\lim_{k \to \infty} \| (w^k;u^k;v^k) \| = \infty$. Hence, the coercivity of $\mathcal{V}_{\rho,\mu} (w,u,\widetilde{z},v)$ implies that $\mathcal{V}_{\rho,\mu} (w^k,u^k,\widetilde{z},v^k) > M$ holds for sufficiently large $k$, which is a contradiction due to $(w^k;u^k;v^k) \in \Omega_M$. Therefore, $\Omega_M$ is a nonempty and bounded set. Then it follows from \cite[Theorem 4.10 (ii)]{mordukhovich2013easy} that (\ref{exact_penalty}) exists a global minimizer. \hfill \Halmos 
\endproof


Generally speaking, the critical point is a nice stationary point for optimality analysis. Particularly, if $(\widehat{w};\widehat{u};\widehat{v})$ is a critical point of (\ref{exact_penalty}), then there exists a Lagrange multiplier $\widehat{z}$ such that
\begin{equation} \label{crit_ori}
	\left\{ \begin{aligned}
		& \nabla f(\widehat{w}) + \mu (\widehat{w} - \widehat{v}) + A^\top \widehat{z} = 0,\\
		& \widehat{z} \in \partial \| \widehat{u}_+ \|_0, \\
		& \widehat{z} = \widetilde{z} + \rho (A\widehat{w} + b - \widehat{u}), \\
		& \widehat{v} = \widehat{w}.
	\end{aligned}  \right.
\end{equation}
Comparing to the critical point, the proximal-type stationary point satisfies a fixed point equation which is a more desirable form for algorithm design. Therefore, we define the following proximal-type stationary point for the optimality analysis of (\ref{exact_penalty}). 


\begin{definition} \label{def_P^*_sta}	A vector $(\widehat{w};\widehat{u};\widehat{v}) \in \mathbb{R}^{p} \times \mathbb{R}^{n} \times \mathbb{R}^{p}$ is called a $\widetilde{\rm P}$-stationary point of (\ref{exact_penalty}) if there exist a vector $\widehat{z} \in \mathbb{R}^n$ and a constant $\alpha >0$ such that  
	\begin{equation} \label{P^*-stat}
		\left\{ \begin{aligned}
			& \nabla f(\widehat{w}) + \mu (\widehat{w} - \widehat{v}) + A^\top \widehat{z} = 0,\\
			& \widehat{u} \in {\rm Prox_{\alpha\lambda \|(\cdot)_+\|_0}} ( \widehat{u} + \alpha \widehat{z} ), \\
			& \widehat{z} = \widetilde{z} + \rho (A\widehat{w} + b - \widehat{u}), \\
			& \widehat{v} = \widehat{w}.
		\end{aligned}  \right.
	\end{equation}
	Moreover, we say $\widehat{z}$ is a $\widetilde{\rm P}$-stationary multiplier of problem (\ref{exact_penalty}).
\end{definition}

Obviously, the difference between (\ref{crit_ori}) and (\ref{P^*-stat}) is the second formula. Actually, the $\widetilde{\rm P}$-stationary point utilizes a positive parameter $\alpha$ to restrict the range of $\widehat{u}$ and $\widehat{z}$. To further clarify the relationship and difference between the critical point and $\widetilde{\rm P}$-stationary point, we give the following property.

\begin{property}\label{rela_cr_p}
	Given $\alpha > 0$, if $\widehat{u}$ and $\widehat{z}$ satisfy $\widehat{u} \in {\rm Prox_{\alpha\lambda \|(\cdot)_+\|_0}} ( \widehat{u} + \alpha \widehat{z} )$, then $\widehat{z} \in \partial \| \widehat{u}_+ \|_0$ holds, while the converse is not true. Therefore, for problem (\ref{exact_penalty}), a $\widetilde{\rm P}$-stationarity with $\alpha > 0$ is strictly stronger than a critical stationarity. 
\end{property}
\proof{Proof.}
Let us consider $\widehat{u} \in {\rm Prox_{\alpha\lambda \|(\cdot)_+\|_0}} ( \widehat{u} + \alpha \widehat{z} )$ with a fixed $\alpha > 0$ from the following three cases.

Case I: If $\widehat{u}_i + \alpha \widehat{z}_i \in (0, \sqrt{2\lambda\alpha})$, then $\widehat{u}_i = {\rm Prox_{\alpha\lambda \|(\cdot)_+\|_0}} ( \widehat{u}_i + \alpha \widehat{z}_i ) = 0$ and hence $\widehat{z}_i \in (0,\sqrt{2\lambda/\alpha})$.

Case II: If $\widehat{u}_i + \alpha \widehat{z}_i \in \{0, \sqrt{2\lambda\alpha}\}$, then $\widehat{u}_i \in {\rm Prox_{\alpha\lambda \|(\cdot)_+\|_0}} ( \widehat{u}_i + \alpha \widehat{z}_i ) = \{ 0, \widehat{u}_i + \alpha \widehat{z}_i \}$. Therefore, $\widehat{u}_i = 0, \widehat{z}_i = 0$ or $\widehat{u}_i = 0, \widehat{z}_i = \sqrt{2\lambda/\alpha}$ or $\widehat{u}_i = \sqrt{2\lambda\alpha}, \widehat{z}_i = 0$.

Case III: If $\widehat{u}_i + \alpha \widehat{z}_i \in (-\infty,0) \cup (\sqrt{2\lambda\alpha}, \infty)$, then $\widehat{u}_i = {\rm Prox_{\alpha\lambda \|(\cdot)_+\|_0}} ( \widehat{u}_i + \alpha \widehat{z}_i ) = \widehat{u}_i + \alpha \widehat{z}_i$ and hence $\widehat{z}_i = 0, \widehat{u}_i \in (-\infty,0) \cup (\sqrt{2\lambda\alpha}, \infty)$.

The above three cases imply that $\widehat{u}$ and $\widehat{z}$ satisfy
\begin{equation}  \label{uz_prox}
	\left\{ \begin{aligned}
		& \widehat{z}_i = 0  , &&\widehat{u}_i \in (-\infty,0] \cup [\sqrt{2\lambda\alpha}, \infty), \\
		& \widehat{u}_i = 0, &&  \widehat{z}_i \in [0,\sqrt{2\lambda/\alpha}],
	\end{aligned}  \right. ~\forall i \in [n].
\end{equation}
From (\ref{sub+0}), $\widehat{z} \in \partial \| \widehat{u}_+ \|_0$ means that $\widehat{u}$ and $\widehat{z}$ satisfy
\begin{equation}\label{uz_diff}
	\left\{ \begin{aligned}
		&\widehat{z}_i=0,  && \widehat{u}_i \neq 0, \\
		&\widehat{u}_i = 0, &&\widehat{z}_i\geq0,
	\end{aligned} \right.~\forall i \in [n].
\end{equation}
It is obvious that (\ref{uz_prox}) yields (\ref{uz_diff}), but the converse does not hold. Taking (\ref{crit_ori}) and (\ref{P^*-stat}) into account, we know that a $\widetilde{\rm P}$-stationarity with $\alpha > 0$ is strictly stronger than a critical stationarity for problem (\ref{exact_penalty}). \hfill\Halmos
\endproof





Now let us conduct optimality analysis on problem (\ref{exact_penalty}).

\begin{theorem} \label{P_sta_optcon}
	The following assertions for the first-order optimality condition of problem (\ref{exact_penalty}) hold.
	
	(i) Given a vector $(\widehat{w};\widehat{u};\widehat{v}) \in \mathbb{R}^{p} \times \mathbb{R}^{n} \times \mathbb{R}^{p}$, if $(\widehat{w};\widehat{u};\widehat{v})$ is a local minimizer, then $(\widehat{w};\widehat{u};\widehat{v})$ is a $\widetilde{\rm P}$-stationary point with a $\widetilde{\rm P}$-stationary multiplier $\widehat{z}$ for any $0 < \alpha < \widehat{\alpha}$, where $ \widehat{\alpha} = \min\{ \widehat{\alpha}_{u}, \widehat{\alpha}_{z} \}$,
	\begin{equation} \notag \label{alpha_hat}
		\widehat{\alpha}_{u} = \left\{  \begin{aligned}
			& +\infty, && \widehat{u}\leq 0, \\
			&\min\limits_{i: \widehat{u}_i > 0}\frac{\widehat{u}_i^2}{ 2\lambda }, && {\rm otherwise},
		\end{aligned} \right. 
		~~~  \widehat{\alpha}_{z} = \left\{  \begin{aligned}
			& +\infty, && \widehat{z}\leq 0, \\
			&\min\limits_{i: \widehat{z}_i > 0}\frac{2\lambda}{ \widehat{z}_i^2 }, && {\rm otherwise}.  
		\end{aligned} \right.  
	\end{equation} 
	
	(ii) If $f$ is convex, then a $\widetilde{\rm P}$-stationary point with $\alpha > 0$ is a local minimizer.
	
	(iii) If $f$ is $\sigma_f$-strongly convex, then a $\widetilde{\rm P}$-stationary point with $\alpha > 1/\sigma$ is the unique global minimizer, where $\sigma$ is defined in (\ref{h_str_con}).
\end{theorem}
\proof{Proof.}
(i) Since $( \widehat{w}; \widehat{u}; \widehat{v} )$ is a local minimizer of (\ref{exact_penalty}), from \citet[Theorem 10.1, Exercise 10.10]{RockWets98}, we have
\begin{equation}
	0 \in \nabla \widetilde{h}(\widehat{w}, \widehat{u},\widehat{v}) + \left(  \begin{array}{c}
		0 \\
		\partial \| \widehat{u}_+ \|_0 \\
		0
	\end{array}  \right) = \left(  \begin{array}{c}
		\nabla_w \widetilde{h}(\widehat{w}, \widehat{u},\widehat{v}) \\
		\nabla_u \widetilde{h}(\widehat{w}, \widehat{u},\widehat{v}) + \partial \| \widehat{u}_+ \|_0 \\
		\nabla_v \widetilde{h}(\widehat{w}, \widehat{u},\widehat{v})
	\end{array}  \right).
\end{equation}
Taking $\widehat{z} = - \nabla_u \widetilde{h}(\widehat{w}, \widehat{u},\widehat{v}) = \widetilde{z} + \rho (A\widehat{w} + b - \widehat{u})  $, then from (\ref{nabla_h3}), we know that (\ref{crit_ori}) holds.
To prove that a vector $(\widehat{w};\widehat{u};\widehat{z})$ satisfying (\ref{crit_ori}) yields (\ref{P^*-stat}) with the $0 < \alpha < \widehat{\alpha}$, it suffices to show that vectors $\widehat{u}$ and $\widehat{z}$ satisfying $\widehat{z} \in \partial \| \widehat{u}_+ \|_0$ yield $\widehat{u} \in {\rm Prox_{\alpha\lambda \|(\cdot)_+\|_0}} ( \widehat{u} + \alpha \widehat{z} )$ with $0 < \alpha < \widehat{\alpha}$.	

For the vector $\widehat{u}$, if there exists $i_0 \in [n]$ with $\widehat{u}_{i_0} >0$, then from $0< \alpha < \widehat{\alpha}$, we have the following estimation
\begin{align} \notag
	\alpha < \widehat{\alpha} \leq \widehat{\alpha}_u = \min\limits_{i: \widehat{u}_i > 0}\frac{\widehat{u}_i^2}{ 2\lambda } \leq  \frac{\widehat{u}_{i_0}^2}{ 2\lambda },
\end{align}
which implies $\widehat{u}_{i_0}>\sqrt{2\lambda\alpha}$. Similarly, if there exists $i_0 \in [n]$ with $\widehat{z}_{i_0} >0$, then from $0< \alpha < \widehat{\alpha}$, we have the following estimation
\begin{align} \notag
	\alpha < \widehat{\alpha} \leq \widehat{\alpha}_z = \min\limits_{i: \widehat{z}_i > 0}\frac{2\lambda}{ \widehat{z}_i^2 } \leq  \frac{2\lambda}{ \widehat{z}_{i_0}^2 },
\end{align}
which implies $\widehat{z}_{i_0} < \sqrt{2\lambda/\alpha}$. Then $\widehat{z} \in \partial \| \widehat{u}_+ \|_0$ yields that for any $0 < \alpha < \widehat{\alpha}$,
\begin{equation} \notag
	\left\{ \begin{aligned}
		&\widehat{z}_i=0,  && \widehat{u}_i \in (-\infty, 0]\cup (\sqrt{2\lambda\alpha}, \infty) , \\
		&\widehat{u}_i = 0, &&\widehat{z}_i \in  [0,\sqrt{2\lambda/\alpha}).
	\end{aligned} \right.
\end{equation}
Now let us consider the following three cases.

Case I: $\widehat{z}_i = 0, \widehat{u}_i \in (-\infty, 0)\cup (\sqrt{2\lambda\alpha}, \infty)$. 
Then by direct calculation, we have $ \widehat{u}_i + \alpha \widehat{z}_i \in (-\infty, 0)\cup (\sqrt{2\lambda\alpha}, \infty)$ and hence $\widehat{u}_i = \widehat{u}_i + \alpha \widehat{z}_i = [ {\rm Prox_{\alpha\lambda \|(\cdot)_+\|_0}} ( \widehat{u} + \alpha \widehat{z} ) ]_i$. 

Case II: $\widehat{z}_i = 0, \widehat{u}_i = 0$. Then $ \widehat{u}_i + \alpha \widehat{z}_i = 0 \in \{0, \sqrt{2\lambda\alpha}\}$; therefore, $\widehat{u}_i \ = 0 \in [ {\rm Prox_{\alpha\lambda \|(\cdot)_+\|_0}} ( \widehat{u} + \alpha \widehat{z} ) ]_i$.

Case III: $\widehat{z}_i \in  (0,\sqrt{2\lambda/\alpha}), \widehat{u}_i = 0$. 
By direct computation, we have $ \widehat{u}_i + \alpha \widehat{z}_i \in (0,\sqrt{2\lambda\alpha})$ and hence $\widehat{u}_i = 0 = [ {\rm Prox_{\alpha\lambda \|(\cdot)_+\|_0}} ( \widehat{u} + \alpha \widehat{z} ) ]_i$. 

Overall, the above three cases imply the conclusion of (i).

(ii) Since $(\widehat{w};\widehat{u};\widehat{v})$ is a $\widetilde{\rm P}$-stationary point of (\ref{exact_penalty}) with $\alpha > 0$, it follows from (\ref{P^*-stat}) and (\ref{nabla_h3}) that
\begin{equation} \label{P^*_stat1}
	\left\{ \begin{aligned}
		& \nabla_w \widetilde{h}( \widehat{w}, \widehat{u},\widehat{v} ) = 0,\\
		& \widehat{u} \in {\rm Prox_{\alpha\lambda \|(\cdot)_+\|_0}} ( \widehat{u} + \alpha \widehat{z} ), \\
		& \widehat{z} = -\nabla_u \widetilde{h}( \widehat{w}, \widehat{u}, \widehat{v} ), \\
		& \nabla_{v} \widetilde{h}(\widehat{w},\widehat{u},\widehat{v}) = 0.
	\end{aligned}  \right.
\end{equation}	
From Property \ref{rela_cr_p}, $\widehat{u} \in {\rm Prox_{\alpha\lambda \|(\cdot)_+\|_0}} ( \widehat{u} + \alpha \widehat{z} )$ with $\alpha > 0$ yields (\ref{uz_prox}).	Denote $\widehat{S}_0 : = \{ i \in [n]: \widehat{u}_i = 0 \}$ and its complementary set $\widehat{S}_0^c : = \{ i \in [n]: \widehat{u}_i \neq 0 \}$,
then (\ref{uz_prox}) yields that
\begin{align} 
	& \widehat{z}_i\widehat{u}_i = 0, ~ \forall i \in [n]. \label{comp_zu} \\
	&0 \leq \widehat{z}_i \leq \sqrt{2\lambda/\alpha},~\forall i \in \widehat{S}_0,~{\rm and}~ \widehat{z}_i = 0,~ \forall i \in \widehat{S}_0^c. \label{z0neq0}
\end{align}
Let us define radius
$\widehat{\delta} := \min\{ \widehat{\delta}_1, \widehat{\delta}_2 \}$, where
\begin{equation} \label{delta}
	\widehat{\delta}_1 := \left\{ \begin{aligned}
		& +\infty, &&{\rm if}~\widehat{z}_{\widehat{S}_0} = 0, \\
		&\sqrt{\frac{\lambda\alpha}{2n}}  , &&{\rm otherwise}, 
	\end{aligned} \right.  
	~~\widehat{\delta}_2 := \left\{ \begin{aligned}
		& +\infty, &&{\rm if}~ \widehat{u} \leq 0, \\
		& \sqrt{2\lambda\alpha}, &&{\rm otherwise}. 
	\end{aligned} \right.
\end{equation}
Then for any $u \in \mathcal{N}(\widehat{u},\delta)$, if $\widehat{u}_i > 0$, we have $\widehat{u}_i > \sqrt{2\lambda\alpha}$ and
\begin{equation} \notag
	u_i = \widehat{u}_i - ( \widehat{u}_i - u_i ) > \widehat{u}_i - \widehat{\delta}_2 \mathop{\geq}\limits^{(\ref{delta})} \sqrt{2\lambda\alpha} - \sqrt{2\lambda\alpha} = 0,
\end{equation}
which implies that 
\begin{align}
	& \{ i \in [n]: \widehat{u}_i > 0 \} \subseteq \{ i \in [n]: u_i > 0 \}. \label{S0+Shat+}
\end{align}

Since $f$ is convex, Property \ref{prop_htilde} (ii) implies that $\widetilde{h}$ is a convex function. Then for any $(w;u;v) \in \mathbb{R}^{p} \times \mathbb{R}^{n} \times \mathbb{R}^{p}$, we have
\begin{align}
	&\widetilde{h}(w,u,v) + \lambda \| u_+ \|_0 - \widetilde{h}(\widehat{w},\widehat{u},\widehat{v}) - \lambda \| \widehat{u}_+ \|_0 \notag \\
	\geq~& \langle \nabla_w \widetilde{h}(\widehat{w},\widehat{u},\widehat{v}), w - \widehat{w}  \rangle +  \langle \nabla_u \widetilde{h}(\widehat{w},\widehat{u},\widehat{v}), u - \widehat{u}  \rangle + \langle \nabla_v \widetilde{h}(\widehat{w},\widehat{u},\widehat{v}), v - \widehat{v}  \rangle + \lambda \| u_+ \|_0 - \lambda \| \widehat{u}_+ \|_0 \notag \\
	\mathop{=}\limits^{(\ref{P^*_stat1})}& \langle  -\widehat{z}, u - \widehat{u} \rangle + \lambda \| u_+ \|_0 - \lambda \| \widehat{u}_+ \|_0 \notag \\
	\mathop{=}\limits^{(\ref{z0neq0})}& \langle  -\widehat{z}_{\widehat{S}_0}, u_{\widehat{S}_0} - \widehat{u}_{\widehat{S}_0} \rangle + \lambda \| u_+ \|_0 - \lambda \| \widehat{u}_+ \|_0 \label{h_con} \\
	\mathop{=}\limits^{(\ref{comp_zu})}& \langle  -\widehat{z}_{\widehat{S}_0}, u_{\widehat{S}_0} \rangle + \lambda \| u_+ \|_0 - \lambda \| \widehat{u}_+ \|_0 . \label{h_con1}
\end{align}

We first consider a trivial case $\widehat{z}_{\widehat{S}_0} = 0$.
For any $( w; u; v ) \in \mathcal{N}( ( \widehat{w}; \widehat{u}; \widehat{v} ), \delta )$, note that $(\ref{S0+Shat+})$ implies $\| u_+ \|_0 \geq \| \widehat{u}_+ \|_0$, then (\ref{h_con1}) directly yields $\widetilde{h}(w,u,v) + \lambda \| u_+ \|_0 \geq \widetilde{h}(\widehat{w},\widehat{u},\widehat{v}) + \lambda \| \widehat{u}_+ \|_0$. For $\widehat{z}_{\widehat{S}_0} \neq 0$, let us consider the following two cases.

Case I: $\| u_+ \|_0 = \| \widehat{u}_+ \|_0$. Notice that (\ref{S0+Shat+}) holds, then for any $\widehat{u}_i \leq 0$, we have $u_i \leq 0$, otherwise $\| u_+ \|_0 \geq \| \widehat{u}_+ \|_0 + 1$. This implies that $u_{\widehat{S}_0} \leq 0$. Together with (\ref{z0neq0}), we have  
\begin{equation}  \label{zu>=0}
	\langle  -\widehat{z}_{\widehat{S}_0}, u_{\widehat{S}_0} \rangle \geq 0.
\end{equation}
Combining (\ref{zu>=0}) with (\ref{h_con1}), we can obtain
\begin{equation} \notag
	\widetilde{h}(w,u,v) + \lambda \| u_+ \|_0 - \widetilde{h}(\widehat{w},\widehat{u}, \widehat{v}) - \lambda \| \widehat{u}_+ \|_0 \geq 0.
\end{equation}

Case II: $\| u_+ \|_0 \neq \| \widehat{u}_+ \|_0$. Taking (\ref{S0+Shat+}) into account, we have $\| u_+ \|_0 \geq \| \widehat{u}_+ \|_0 + 1$. It follows from (\ref{h_con}) that
\begin{align}
	\widetilde{h}(w,u,v) + \lambda \| u_+ \|_0 - (\widetilde{h}(\widehat{w},\widehat{u},\widehat{v}) + \lambda \| \widehat{u}_+ \|_0) \geq  \lambda -  \|\widehat{z}_{\widehat{S}_0}\|\| \| u - \widehat{u} \| \mathop{>}\limits^{(\ref{z0neq0})}  \lambda - \widehat{\delta}\sqrt{2n\lambda/\alpha}  
	\mathop{\geq}\limits^{(\ref{delta})}  0. \notag
\end{align}
In conclusion, the above two cases show that $( \widehat{w}; \widehat{u};\widehat{v} )$ is a local minimizer of (\ref{exact_penalty}).	

(iii) Since $(\widehat{w};\widehat{u};\widehat{v})$ is a $\widetilde{\rm P}$-stationary point of (\ref{exact_penalty}) with $\alpha > 0$, then (\ref{P^*_stat1}) holds. By the definition of ${\rm Prox_{\alpha\lambda \|(\cdot)_+\|_0}} ( \cdot )$, for any $u$ and $z \in \mathbb{R}^n$, we have
\begin{equation} \notag
	\frac{1}{2\alpha} \| \widehat{u} - ( \widehat{u} + \alpha \widehat{z} ) \|^2 + \lambda \| \widehat{u}_+ \|_0 \leq \frac{1}{2\alpha} \| u - ( \widehat{u} + \alpha \widehat{z} ) \|^2 + \lambda \| u_+ \|_0.
\end{equation}
After some simple algebraic manipulation, we can obtain
\begin{equation}  \label{diff_+0}
	\lambda \| u_+ \|_0 - \lambda \| \widehat{u}_+ \|_0 \geq \langle \widehat{z}, u - \widehat{u} \rangle - \frac{1}{2\alpha} \| u - \widehat{u} \|^2.
\end{equation}
From the $\sigma_f$-strongly convexity of $f$ and Property \ref{prop_htilde} (ii), we know that $\widetilde{h}$ is $\sigma$-strongly convex. 
For any $(w;u;v) \in \mathbb{R}^{p} \times \mathbb{R}^{n} \times \mathbb{R}^{p}$ and $(w;u;v) \neq (\widehat{w};\widehat{u};\widehat{v})$, we have the following estimate
\begin{align}
	&\widetilde{h}(w,u,v) + \lambda \| u_+ \|_0 - \widetilde{h}(\widehat{w},\widehat{u},\widehat{v}) - \lambda \| \widehat{u}_+ \|_0 \notag\\
	\mathop{\geq}\limits^{(\ref{h_str_con})}& \langle \nabla_w \widetilde{h}(\widehat{w},\widehat{u},\widehat{v}), w - \widehat{w}  \rangle + \langle \nabla_u \widetilde{h}(\widehat{w},\widehat{u},\widehat{v}), u - \widehat{u}  \rangle + \langle \nabla_v \widetilde{h}(\widehat{w},\widehat{u},\widehat{v}), v - \widehat{v}  \rangle \notag \\
	& + \frac{\sigma}{2} (\| w - \widehat{w} \|^2 + \| u - \widehat{u} \|^2 + \| v - \widehat{v} \|^2 )+ \lambda \| u_+ \|_0 - \lambda \| \widehat{u}_+ \|_0 \notag \\
	\mathop{=}\limits^{(\ref{P^*_stat1})}& -\langle \widehat{z}, u - \widehat{u} \rangle+ \frac{\sigma}{2} (\| w - \widehat{w} \|^2 + \| u - \widehat{u} \|^2 + \| v - \widehat{v} \|^2  ) + \lambda \| u_+ \|_0 - \lambda \| \widehat{u}_+ \|_0 \notag \\
	\mathop{\geq}\limits^{(\ref{diff_+0})}& \frac{\sigma}{2} \| w - \widehat{w} \|^2 +  \frac{1}{2}(\sigma - \frac{1}{\alpha}) \| u - \widehat{u} \|^2  + \frac{\sigma}{2} \| v - \widehat{v} \|^2 > 0. \notag
\end{align}
The last inequality follows from $\alpha > 1/\sigma$ and $(w;u;v) \neq (\widehat{w};\widehat{u};\widehat{v})$. This shows that $(\widehat{w};\widehat{u};\widehat{v})$ is the unique global minimizer of (\ref{exact_penalty}). \hfill\Halmos
\endproof

\section{Lyapunov (Strongly) Exact Penalty.} \label{lya_ep}
In this section, we will show the (strongly) exact penalty properties of Lyapunov function for 0/1-COP. The Lyapunov exact penalty means that  a solution (minimizer or proximal-type stationary point) of (\ref{eq2}) yields a solution of the Lyapunov function with an appropriately selected Lagrange multiplier when $\rho$ and $\mu$ are sufficiently large, meanwhile the Lyapunov strongly exact penalty means that the converse inclusion holds, which is a more difficult but desirable result. Let us begin this section with the optimality condition of (\ref{eq2}). Zhou et al. \cite{zhou2021quadratic} introduced the concept of P-stationary point for establishing the first-order optimality condition of (\ref{eq2}).

\begin{definition} \label{def_P_sta}
	A vector $(w^*;u^*) \in \mathbb{R}^{p} \times \mathbb{R}^{n}$ is called a {\rm P}-stationary point of (\ref{eq2}) if there exist a vector $z^* \in \mathbb{R}^n$ and a constant $\alpha >0$ such that  
	\begin{equation} \label{P-stat}
		\left\{ \begin{aligned}
			& \nabla f(w^*) + A^\top z^* = 0,\\
			& u^* \in {\rm Prox_{\alpha\lambda \|(\cdot)_+\|_0}} ( u^* + \alpha z^* ), \\
			& Aw^* + b - u^* = 0.
		\end{aligned}  \right.
	\end{equation}
	Moreover, we say $z^*$ is a {\rm P}-stationary multiplier and $(w^*;u^*;z^*)$ is a {\rm P}-stationary triplet of (\ref{eq2}).
\end{definition}

Given a feasible couple $(w^*;u^*)$ of (\ref{eq2}), we define $S^*_0 := \{ i \in [n] : u^*_i = 0 \}$. The following theorem proposed in \citet{zhou2021quadratic} shows the relationship between a P-stationary point and a local (global) minimizer of (\ref{eq2}).

\begin{lemma}\textbf{(\citet[Theorem 3.3]{zhou2021quadratic})} \label{opt_orig}
	The following assertions for the first-order optimality condition for problem (\ref{eq2}) hold.
	
	(i) Given a vector $(w^*;u^*) \in \mathbb{R}^{p} \times \mathbb{R}^{n}$, if $A_{S^*_0}$ is full row rank and $(w^*,u^*)$ is a local minimizer, then $(w^*,u^*)$ is a {\rm P}-stationary point and with a {\rm P}-stationary multiplier $z^*$ for any $0 < \alpha < \alpha^*$, where $\alpha^* := \min\{ \alpha^*_{u}, \alpha^*_{z} \}$,
	\begin{equation} \notag
		\alpha^*_{u} := \left\{  
		\begin{aligned}
			& +\infty,  &&u^*\leq 0, \\
			&\min\{\frac{{u^*_i}^2}{ 2\lambda }: u_i^* > 0 \}, &&{\rm otherwise},
		\end{aligned} \right. , ~
		\alpha^*_{z} = \left\{  \begin{aligned}
			& +\infty, && z^*\leq 0, \\
			&\min\{\frac{2\lambda}{ {z^*_i}^2 } : z^*_i > 0 \}, && {\rm otherwise}. 
		\end{aligned} \right.  
	\end{equation} 
	
	(ii) If $f$ is a convex function, then a {\rm P}-stationary point with $\alpha > 0$ is a local minimizer.
	
	(iii) If $f$ is a $\sigma_f$-strongly convex function, then a {\rm P}-stationary point with $\alpha > \|A\|^2/\sigma_f$ is the unique global minimizer.
\end{lemma}

We stress that Lemma \ref{opt_orig} (iii) is slightly different with the Theorem 3.3 (iii) in \citet{zhou2021quadratic}. Under the assumption that $f$ is strongly convex, Zhou et al. proved that a P-stationary point with $\alpha \geq \|A\|^2/\sigma_f$ is a global minimizer of (\ref{eq2}). Actually, taking the similar procedure to the proof of Theorem 3.3 (iii) in \citet{zhou2021quadratic}, we can obtain that a P-stationary point with $\alpha > \|A\|^2/\sigma_f$ is the unique global minimizer of (\ref{eq2}).


Given a solution $(w^*;u^*)$ of (\ref{eq2}), Lemma \ref{opt_orig} shows that there exists a {\rm P}-stationary multiplier $z^*$. We take $\widetilde{z} = z^*$ and reformulate (\ref{exact_penalty}) in the following problem
\begin{equation} \label{exact_pen_eq}
	\min_{(w;u;v) \in \mathbb{R}^{p} \times \mathbb{R}^{n} \times \mathbb{R}^{p}} \mathcal{V}_{\rho,\mu} ( w,u,z^*,v).
\end{equation}
Similarly, we can 
define a ${\rm P}^*$-stationary point and a ${\rm P}^*$-stationary multiplier for problem (\ref{exact_pen_eq}). Therefore, (\ref{exact_pen_eq}) has the same optmality condition as (\ref{P^*-stat}).

In the subsequent theorem, we will show that a solution of (\ref{eq2}) yields a solution of (\ref{exact_pen_eq}) under some conditions, namely $\mathcal{V}_{\rho,\mu} ( w,u,z^*,v )$ is an exact penalty function of 0/1-COP.


\begin{theorem}[Exact Penalty Theorem] \label{exact_penalty_relation}The following assertions hold for problems (\ref{eq2}) and (\ref{exact_pen_eq}).
	
	(i) Given a {\rm P}-stationary point $(w^*;u^*)$ of (\ref{eq2}) with $\alpha > 0$ and ${\rm P}$-stationary multiplier $z^*$ , then $(w^*;u^*;w^*)$ is a ${\rm P}^*$-stationary point with $\alpha > 0$ of (\ref{exact_pen_eq}) for any $\mu > 0$ and $\rho > 0$.
	
	(ii) Given a local minimizer $(w^*;u^*)$ of (\ref{eq2}), if $A_{S_0^*}$ is full row rank and $f$ is convex, then there exists a ${\rm P}$-stationary multiplier $z^*$ of (\ref{eq2}) such that $(w^*;u^*;w^*)$ is a local minimizer of (\ref{exact_pen_eq}) for any $\mu > 0$ and $\rho > 0$.
	
	(iii) Given a {\rm P}-stationary point $(w^*;u^*)$ of (\ref{eq2}) with $\alpha > (\| A \|^2 + 2)/ \sigma_f$ and ${\rm P}$-stationary multiplier $z^*$, if $f$ is $\sigma_f$-strongly convex and the parameters $\mu$ and $\rho$ are set as
	\begin{equation} \label{para_exact_penalty}
		\mu > \frac{\sigma_f}{\alpha\sigma_f - \| A \|^2 - 2},~\rho > \frac{\mu\sigma_f}{\mu( \alpha\sigma_f - \| A \|^2 - 2 ) - \sigma_f},
	\end{equation}
	then $(w^*;u^*;w^*)$ is the unique global minimizer of (\ref{exact_pen_eq}).
\end{theorem} 

\proof{Proof.}
(i)	Since $(w^*;u^*)$ is a {\rm P}-stationary point of (\ref{eq2}) with {\rm P}-stationary multiplier $z^*$, then from (\ref{P-stat}), we have 
\begin{equation} \notag 
	\left\{ \begin{aligned}
		& \nabla f(w^*) + \mu (w^* - v^*) + A^\top z^* = 0,\\
		& u^* \in {\rm Prox_{\alpha\lambda \|(\cdot)_+\|_0}} ( u^* + \alpha z^*), \\
		& z^* = z^* + \rho (Aw^* + b - u^*), \\
		&v^* = w^*.
	\end{aligned}  \right.
\end{equation}
This implies that $(w^*;u^*;w^*)$ is a ${\rm P}^*$-stationary point of (\ref{exact_pen_eq}).

(ii) Since $(w^*;u^*)$ is a local minimizer of (\ref{eq2}) and $A_{S_0^*}$ is full row rank, it follows from Lemma \ref{opt_orig} (i) that there exists a P-stationary multiplier $z^*$ such that (\ref{P-stat}) holds with $0 < \alpha < \widehat{\alpha}$. Then assertion (i) yields that $(w^*;u^*;w^*)$ is a ${\rm P}^*$-stationary point of (\ref{exact_pen_eq}) with $0 < \alpha < \widehat{\alpha}$. Theorem \ref{P_sta_optcon} (ii) directly implies that $(w^*;u^*;w^*)$ is a local minimizer of (\ref{exact_pen_eq}).

(iii) Assertion (i) implies that  $(w^*;u^*;w^*)$ is a ${ \rm P}^*$-stationary point of (\ref{exact_pen_eq}) with a positive constant $\alpha > (\| A \|^2 + 2)/ \sigma_f$. Since parameters $\mu,\rho$ are chosen as (\ref{para_exact_penalty}), direct computation implies that 
\begin{equation} \notag
	\alpha > \frac{\sigma_f(\mu + \rho ) + \rho\mu \| A \|^2 + 2\rho \mu}{\sigma_f \rho \mu} = \frac{1}{\sigma}.
\end{equation}
Theorem \ref{P_sta_optcon} (iii) directly yields that $(w^*;u^*;w^*)$ is the unique global minimizer of (\ref{exact_pen_eq}). 
\hfill	\Halmos 
\endproof

Theorem \ref{exact_penalty_relation} can be seen as an extension of the classical result for the smooth nonlinear programming (see e.g. \citet[Theorem 17.5]{nocedal2006numerical}). In practice, the converse inclusions of Theorem \ref{exact_penalty_relation} is more desirable. Therefore, we give the following theorem.

\begin{theorem}[Strongly Exact Penalty Theorem] \label{exa_pen2}
	Given a ${\rm \widetilde{P}}$-stationary point $(\widehat{w};\widehat{u};\widehat{v})$ of (\ref{exact_penalty}) with ${\rm \widetilde{P}}$-stationary multiplier $\widehat{z}$ and constant $\alpha > 0$, 
	then the following assertions hold.
	
	(i) $(\widehat{w};\widehat{u})$ is a ${ \rm P }$-stationary point of (\ref{eq2}) with ${ \rm P }$-stationary multiplier $\widehat{z}$ and constant $\alpha > 0$ if and only if $ \widehat{z} = \widetilde{z}$.
	
	(ii) If $A$ is full row rank, then $\widehat{z} = -(AA^\top)^{-1} A \nabla f(\widehat{w})$. Furthermore, $(\widehat{w};\widehat{u})$ is a ${ \rm P }$-stationary point of (\ref{eq2}) with ${ \rm P }$-stationary multiplier $\widehat{z}$ and $\alpha > 0$ if and only if $\widehat{z}  = -(AA^\top)^{-1} A \nabla f(\widehat{w})= \widetilde{z}$.
\end{theorem} 

\proof{Proof.}
(i) `` $\Rightarrow$ " Since $(\widehat{w};\widehat{u};\widehat{v})$ is a ${\rm \widetilde{P}}$-stationary point of (\ref{exact_penalty}) with ${\rm \widetilde{P}}$-stationary multiplier $\widehat{z}$, then $\widehat{z} - \widetilde{z} = \rho (A \widehat{w} + b - \widehat{u})$ by (\ref{P^*-stat}). Moreover, $(\widehat{w};\widehat{u})$ is a ${ \rm P }$-stationary point of (\ref{eq2}) implies that $A \widehat{w} + b - \widehat{u} = 0$. Hence, $\widehat{z} = \widetilde{z}$ holds.

``$\Leftarrow$" Since $(\widehat{w};\widehat{u};\widehat{v})$ is a ${\rm \widetilde{P}}$-stationary point of (\ref{exact_penalty}) with ${\rm \widetilde{P}}$-stationary multiplier $\widehat{z}$ and constant $\alpha > 0$, and $\widehat{z} = \widetilde{z}$ holds, these facts imply that
\begin{equation} 
	\left\{ \begin{aligned}
		& \nabla f(\widehat{w}) + A^\top \widehat{z} = 0,\\
		& \widehat{u} \in {\rm Prox_{\alpha\lambda \|(\cdot)_+\|_0}} ( \widehat{u} + \alpha \widehat{z} ), \\
		& A\widehat{w} + b - \widehat{u} = 0,
	\end{aligned}  \right.
\end{equation}
and the conclusion of assertion (i) holds.  

(ii) Since $(\widehat{w};\widehat{u};\widehat{v})$ is a ${\rm \widetilde{P}}$-stationary point of (\ref{exact_penalty}) with ${\rm \widetilde{P}}$-stationary multiplier $\widehat{z}$, then $\nabla f(\widehat{w}) + A^\top \widehat{z} = 0$ holds by (\ref{P^*-stat}). From the full row rank of $A$, we know that $\widehat{z} = - (AA^\top)^{-1}A \nabla f(\widehat{w})$. Then taking the similar procedure of assertion (i), we can arrive at the conclusion of assertion (ii).
\hfill	\Halmos 
\endproof

Theorem \ref{exa_pen2} shows that by an appropriated choice of $\widetilde{z}$, a ${\rm \widetilde{P}}$-stationary point of (\ref{exact_penalty}) yields a ${ \rm P }$-stationary point of (\ref{eq2}). However, Theorem \ref{exa_pen2} also indicates that the appropriated $\widetilde{z}$ could not be identified in advance. In other words, we could not obtain the solution of 0/1-COP by minimizing a Lyapunov function $\mathcal{V}_{\rho,\mu} (w,u,\widetilde{z},v)$ with arbitrarily selected $\widetilde{z}$. As an alternative, we consider designing an algorithm for generating a sequence of $\{ z^k \}_{k \in \mathbb{N}}$ to approximate the optimal Lagrange multiplier of (\ref{eq2}).



\section{Inexact Augment Lagrangian Method for  0/1-COP.} \label{ialm}
In this section, we focus on designing an inexact augmented Lagrangian method (IALM) for solving 0/1-COP. Denote $(w^k;u^k;z^k;v^k)$ as the $k$th iteration point. Let us define function
\begin{equation} \notag
	h_k(w,u):= f(w) + \frac{\mu}{2} \| w - v^{k} \|^2 +  \langle z^k, Aw + b - u \rangle + \frac{\rho}{2} \| Aw + b - u \|^2.
\end{equation}
The framework of IALM is as follows.
\begin{algorithm}[H]
	\caption{IALM: an inexact augmented Lagrangian method} \label{IALM}
	\begin{algorithmic}
		\STATE{Initialization: Given a positive sequence $\{\epsilon_k\}_{k \in \mathbb{N}}$, fix $\mu > 0$, $\rho > 0$, $\alpha = 1/\rho$, start with any $(w^0;u^0;z^0;v^0) \in \mathbb{R}^{p} \times \mathbb{R}^{n} \times \mathbb{R}^{n} \times \mathbb{R}^{p}$.}
		\FOR{$k = 0,1, \cdots$}
		\STATE{\textbf{Primal step} Use 0/1-BALM (see Algorithm \ref{pam}) to compute $(w^{k+1};u^{k+1})$ as an approximate solution of
			\begin{equation} \label{eq2.1}
				\mathop{\min}\limits_{w \in \mathbb{R}^p, u \in \mathbb{R}^n} \mathcal{V}_{\rho,\mu} ( w,u,z^k,v^{k} )
			\end{equation}
			satisfying 
			\begin{subequations} \label{stopping1}
				\begin{align} 
					&\mathcal{V}_{\rho,\mu} (w^{k+1},u^{k+1},z^k,v^k) \leq \mathcal{V}_{\rho,\mu} (w^{k},u^{k},z^k,v^k), \label{des_auglag} \\
					&\max \{ {\rm dist}( u^{k+1}, {\rm Prox_{\alpha\lambda \|(\cdot)_+\|_0}} ( u^{k+1} - \alpha \nabla_u h_k ( w^{k+1},u^{k+1} ) ) ), \| \nabla_w h_k(w^{k+1},u^{k+1}) \|   \} \leq \epsilon_k. \label{crit1} 
				\end{align}
		\end{subequations}}
		\STATE{\textbf{Multiplier step}\begin{equation} \label{eq2.2}
				z^{k+1} = z^k + \rho(Aw^{k+1} + b - u^{k+1}).
		\end{equation}}
		\STATE{\textbf{Lyapunov step}  \begin{align}
				v^{k+1} = w^{k+1}. \notag
		\end{align}}
		\ENDFOR
	\end{algorithmic}
\end{algorithm}







To proceed, we give the following assumption on $f$.

\begin{assumption}\label{bbs} The function $f$ is $l_f$-smooth.
\end{assumption}

Comparing to $\widetilde{h}$, the function $h_k$ with fixed Lagrange multiplier $z^k$ and Lyapunov variable $v^k$ has better properties, which are summarized as follows.
\begin{property}\label{prop_hstar}The function $h_k$ has the following properties:
	
	(i) The gradient of $h_k$ can be represented as
	\begin{equation} 
		\begin{aligned} \label{nabla_h}
			\nabla h_k(w,u) &= \left(  \begin{array}{c}
				\nabla_w h_k(w,u) \\
				\nabla_u h_k(w,u)
			\end{array}  \right)  = \left(  \begin{array}{c}
				\nabla f(w) + \mu(w - v^{k}) + \rho A^\top \big( Aw + b - u  + z^k / \rho \big) \\
				\rho(  u - Aw - b  - z^k / \rho )
			\end{array}  \right).
		\end{aligned}
	\end{equation}
	
	
	(ii) If Assumption \ref{bbs} holds, then $ h_k(w,u)$ is $l_h$-smooth, where $l_h := \sqrt{2} \max \{ \mu + l_f + \rho \| A \| ( \|A\| + 1 ), \rho ( \|A\| + 1 ) \}$. 
	
	(iii) If $f$ is convex, then $h_k$ is strongly convex.
	
\end{property}

\proof{Proof.}
Since (i) is easy to be verified, we only prove (ii) and (iii).

(ii) For any $(w;u)$ and $(\overline{w};\overline{u}) \in \mathbb{R}^{p}\times \mathbb{R}^{n}$, we have the following estimation from (\ref{nabla_h})
\begin{align}
	\| \nabla h_k(w,u) - \nabla h_k(\overline{w},\overline{u}) \| =& \left\| \left(  \begin{array}{c}
		\nabla f(w)- \nabla f(\overline{w}) + \mu(w - \overline{w}) + \rho A^\top A (w - \overline{w}) - \rho A^\top ( u - \overline{u} ) \\
		\rho(u - \overline{u}) - \rho A (w - \overline{w})
	\end{array}  \right)
	\right\| \notag \\
	\leq & \| \nabla f(w)- \nabla f(\overline{w}) \| + (\mu + \rho \| A \|^2 + \rho \| A \|) \| w - \overline{w} \| + ( \rho + \rho \| A \| ) \| u - \overline{u} \| \notag \\
	\leq & \max \{ \mu + l_f + \rho \| A \| ( \|A\| + 1 ), \rho ( \|A\| + 1 ) \} ( \| w - \overline{w} \| + \| u - \overline{u} \| ) \notag \\
	\leq& l_h \| (w - \overline{w}; u - \overline{u}) \| \notag,
\end{align}
where the last inequality follows from $t_1 + t_2 \leq \sqrt{2 t_1^2 + 2 t_2^2 } $ for any real numbers $t_1$ and $t_2$.

(iii) First, let us define 
\begin{equation} \label{rep_varphi}
	\phi_k( w,u ): = \frac{\mu}{2} \| w - v^{k} \|^2 +  \langle z^k, Aw + b - u \rangle + \frac{\rho}{2} \| Aw + b - u \|^2. \notag
\end{equation}
Then we have
\begin{align}
	\nabla^2 \phi_k(w,u)  &= \left[ \begin{array}{cc}
		\mu  I + \rho A^\top A   &-\rho A^\top \\
		-\rho A  &\rho I
	\end{array} \right]  = \left[ \begin{array}{cc}
		I   & -  A^\top \\
		0  &I
	\end{array} \right] \left[ \begin{array}{cc}
		\mu  I  &0 \\
		0  &\rho I
	\end{array} \right]  \left[ \begin{array}{cc}
		I   & 0 \\
		-  A  &I
	\end{array} \right] \succ 0,  \label{h_sig_strcon} \notag
\end{align}
which means that $\phi_k$ is strongly convex. Since $f$ is convex and $h_k (w,u) = f(w) + \phi_k (w,u)$, we can conclude that $h_k$ is strongly convex. \hfill \Halmos
\endproof

The problem (\ref{eq2.1}) can be reformulated as 
\begin{align}
	\mathop{\min}\limits_{w \in \mathbb{R}^p, u \in \mathbb{R}^n} \mathcal{V}_{\rho,\mu} ( w,u,z^k,v^{k} ) = h_k(w,u) + \lambda \| u_+ \|_0. \notag
\end{align}
For this problem, similar to (\ref{P^*-stat}), we can define a ${\rm P}^k$-stationary point $(\widehat{w};\widehat{u})$ and a ${\rm P}^k$-stationary multiplier $\widehat{z}$ with $\alpha > 0$ satisfying the following formula
\begin{equation} \label{P^k-stat}
	\left\{ \begin{aligned}
		& \nabla f(\widehat{w}) + \mu (\widehat{w} - v^{k}) + A^\top \widehat{z} = 0,\\
		& \widehat{u} \in {\rm Prox_{\alpha\lambda \|(\cdot)_+\|_0}} ( \widehat{u} + \alpha \widehat{z} ), \\
		& \widehat{z} = z^k + \rho (A\widehat{w} + b - \widehat{u}). \\
	\end{aligned}  \right.
\end{equation}

For solving the primal step, we design a zero-one Bregman alternating linearized minimization method (0/1-BALM). Particularly, in the $(k+1)$th iteration of IALM, 0/1-BALM generates a sequence $\{ (w^{k,j};u^{k,j}) \}_{j \in \mathbb{N}}$ with starting point $(w^{k,0};u^{k,0}) := (w^k;u^k)$. Moreover, we define function
\begin{align} \notag
	\varphi_{k,j}(w):= h_k(w,u^{k,j+1}) + \mathcal{D}_{\psi_{k,j}}(w, w^{k,j}),
\end{align}
where
$\mathcal{D}_{\psi_{k,j}}$ is a Bregman distance with respect to function $\psi_{k,j}$, which means that
\begin{align} \notag
	\mathcal{D}_{\psi_{k,j}} (w, w^{k,j}) := \psi_{k,j}(w) - \psi_{k,j}(w^{k,j}) - \langle \nabla \psi_{k,j}(w^{k,j}), w - w^{k,j} \rangle .
\end{align}
The function $\psi_{k,j}:\mathbb{R}^p \to \mathbb{R}$ is defined as   
\begin{align} \notag
	\psi_{k,j}(w) := \frac{1}{2} \| w \|^2_{P_{k,j}} - g(w),
\end{align}
where $P_{k,j} := tI - Q_{k,j}$, $t > 0$, $Q_{k,j} \succeq 0$ and $g:\mathbb{R}^p \to \mathbb{R}$ is a $l_g$-smooth function. Then the Bregman distance can be further simplified as
\begin{align} \label{breg_sim}
	\mathcal{D}_{\psi_{k,j}} (w, w^{k,j}) = \frac{1}{2} \| w - w^{k,j} \|^2_{P_{k,j}} -( g(w) - g(w^{k,j}) - \langle \nabla g(w^{k,j}), w - w^{k,j} \rangle ).
\end{align}
The Bregman distance is a common strategy to generate a simpler iteration (see e.g. \citet{wang2018con}). In Remark \ref{re-palm}, we will show that an appropriate selection of $Q_{k,j}$ and $g$ simplifies the w-step and improves the performance of our algorithm.  


The framework of 0/1-BALM for solving the subproblem (\ref{eq2.1}) is as follows.

\begin{algorithm}[H] 
	\caption{0/1-BALM: zero-one Bregman alternating linearized minimization algorithm} \label{pam}
	\begin{algorithmic}
		\STATE{Initialization: Take $\alpha = 1/\rho$ and start with $(w^{k,0};u^{k,0}) = (w^k;u^k)$.}
		\FOR{$j=0,1,\cdots$}
		\STATE{\textbf{u-step:} } Compute
		\begin{equation} \label{u-step}
			u^{k,j+1} \in {\rm Prox}_{ \alpha\lambda \| (\cdot)_+ \|_0} ( u^{k,j} - \alpha \nabla_u h_k( w^{k,j}, u^{k,j} ) ).
		\end{equation}
		\STATE{\textbf{w-step:} } Solve the following problem
		\begin{equation} \label{w-step}
			w^{k,j+1} \in \mathop{\arg\,\min}\limits_{w \in \mathbb{R}^p} \varphi_{k,j}(w).
		\end{equation}
		\ENDFOR
	\end{algorithmic}
\end{algorithm}

\begin{remark} \label{re-palm}
	Some details of 0/1-BALM are given as follows.
	\begin{itemize}
		\item[(i)] The iterative formulas (\ref{u-step}) and (\ref{w-step}) can be simplified further. For the u-step (\ref{u-step}), substituting $\alpha = 1/\rho$ and $\nabla_u h_k( w^{k,j}, u^{k,j} ) = - z^k - \rho( Aw^{k,j} + b - u^{k,j} )$ into (\ref{u-step}), we can obtain
		\begin{equation} \label{u-step-sim1}
			u^{k,j+1} \in {\rm Prox}_{ \alpha\lambda \| (\cdot)_+ \|_0} ( Aw^{k,j} + b + z^k/ \rho ).
		\end{equation}
		Denote ${T_{k,j}} := \{ i \in [n]:  0< (Aw^{k,j} + b + z^k/\rho )_i  <  \sqrt{2\lambda\alpha} \}$. Then by (\ref{prox0+}) and (\ref{u-step-sim1}), we could take
		\begin{align}  \label{u_ind}
			u^{k,j+1} = \left( \begin{array}{c}
				u^{k,j+1}_{T_{k,j}} \\
				u^{k,j+1}_{T_{k,j}^c}
			\end{array} \right) =  \left( \begin{array}{c}
				0 \\
				(Aw^{k,j} + b + z^k /\rho)_{T_{k,j}^c}
			\end{array} \right).
		\end{align}
		For the w-step (\ref{w-step}), from the optimality condition, we can obtain
		\begin{align}
			0 =& \nabla \varphi_{k,j}( w^{k,j+1} ) \notag \\
			=&	\nabla f(w^{k,j+1}) + \mu ( w^{k,j+1} - v^{k} ) + \rho A^\top ( Aw^{k,j+1} + b - u^{k,j+1} + \frac{1}{\rho}z^k ) + P_{k,j}( w^{k,j+1}  - w^{k,j}) \notag \\
			& - (\nabla g(w^{k,j+1}) - \nabla g(w^{k,j}) ). \label{opt_phi}
		\end{align}
		Here we give two cases for choosing $Q_{k,j}$ and $g$ to simplify the computation of $w^{k,j+1}$.
		
		Case I: $Q_{k,j} = \rho A^\top A$ and $g(w) = f(w)$. Substituting (\ref{u_ind}), $Q_{k,j}$ and $g$ into (\ref{opt_phi}), we can obtain
		\begin{align} \label{sim_u_step}
			w^{k,j+1} = \frac{\mu}{\mu+t}v^{k} + \frac{t}{t+\mu} w^{k,j} - \frac{1}{\mu+t}( \nabla f(w^{k,j}) + \rho A^\top_{T_{k,j}} ( Aw^{k,j} + b + \frac{1}{\rho}z^k )_{T_{k,j}} ).
		\end{align}
		In this case, the iterative formula of our 0/1-BALM is similar to that of \citet{bolte2014proximal,bolte2018nonconvex,boct2020proximal} and we will illustrate difference and relationship between 0/1-BALM and other methods in the subsequent (ii) and (iii).
		
		Case II: $Q_{k,j} = \rho A^\top_{T_{k,j}^c} A_{T_{k,j}^c}$ and $g(w) = 0$. It should be noted that this setting is particularly suitable for a quadratic function $f$ with a nonsingular and simple Hessian. Without loss of generality, we consider $f(w) = \| w \|^2/2$. Then, from (\ref{opt_phi}), we can obtain
		\begin{align} \notag
			w^{k,j+1} = [ (1+t+\mu)I + \rho A^\top_{T_{k,j}} A_{T_{k,j}}  ]^{-1} ( \mu v^{k} + tw^{k,j} - \rho A^\top_{T_{k,j}} (b+\frac{1}{\rho}z^k)_{T_{k,j}} ).  
		\end{align}
		Taking $(w^{k,0};u^{k,0}) = (w^k;u^k)$ into account, 0/1-ADMM proposed by \citet{wang2021support} can be seen as a degenerate case of our IALM equipped with 0/1-BALM performing only one iteration in the setting of $\mu = t = 0$.
		Although the numerical experiments demonstrate the superior of 0/1-ADMM on computational speed and accuracy for solving SVM, \citet{wang2021support} did not give a global convergence analysis. In Theorems \ref{sub-fin-ter}, \ref{sub_con} and \ref{whole_con}, we will show that our IALM equipped with a 0/1-BALM 
		could ensure the global convergence under suitable conditions.
		
		
		
		(ii) Our 0/1-BALM has lower computational complexity compared with the PALM proposed by \citet{bolte2014proximal}. Specifically, the PALM adopted $\alpha < 1/\rho$ to ensure the sufficient descent of $\mathcal{V}_{\rho,\mu}$ with respect to $w$ and $u$, while we adopt $\alpha = 1/\rho$ in each iteration, which will simplify the computation for w-step. Because the $Aw^{k,j} + b + z^k/\rho$ term has been computed in the u-step, then for the w-step, the computational complexity of our 0/1-BALM is $O(| T_{k,j} | p)$ ($ | T_{k,j} | \leq n $), while the computational complexity of PALM is $O(np)$. This will significantly simplify the computation since a small-scale $T_{k,j}$ could be selected in practice. 
		
		(iii) Our IALM equipped with 0/1-BALM degenerates into the PL-ADMM proposed in \citet{bolte2018nonconvex,boct2020proximal} under specific conditions and has a weaker assumption on $A$ compared with the PL-ADMM. Let us recall the framework of PL-ADMM:
		\begin{align}
			& u^{k+1} \in {\rm Prox}_{ \alpha\lambda \| (\cdot)_+ \|_0} ( Aw^{k} + b + z^k/ \rho ), \notag \\ 
			& w^{k+1} = w^k - \frac{1}{\mu} ( \nabla f(w^k) + \rho A^\top ( A w^k + b - u^{k+1} + z^k / \rho ) ), \notag \\
			&z^{k+1} = z^k + \rho ( Aw^{k+1} + b - u^{k+1} ). \notag
		\end{align}
		Denote $T_k := \{ i \in [n]: 0 < (Aw^{k} + b + z^k/\rho )_i  <  \sqrt{2\lambda\alpha} \}$. Similar to (ii), taking advantage of ${\rm Prox}_{ \alpha\lambda \| (\cdot)_+ \|_0}$, the above framework can be simplified as
		\begin{align}
			& u^{k+1} \in {\rm Prox}_{ \alpha\lambda \| (\cdot)_+ \|_0} ( Aw^{k} + b + z^k/ \rho ), \label{pl-u-step} \\ 
			& w^{k+1} = w^k - \frac{1}{\mu} ( \nabla f(w^k) + \rho A^\top_{T_k} ( A w^k + b - u^{k+1} + z^k / \rho )_{T_k} ), \label{pl-w-step} \\
			&z^{k+1} = z^k + \rho ( Aw^{k+1} + b - u^{k+1} ). \notag
		\end{align}
		It is not hard to see that if we replace the $w^{k,j}$, $u^{k,j}$ and $T_{k,j}$ in (\ref{u-step-sim1}) and (\ref{sim_u_step}) with $w^k$, $u^k$ and $T_k$, then they are similar to (\ref{pl-u-step}) and (\ref{pl-w-step}) respectively. Moreover, note that for the $k$th iteration of our IALM, we take $(w^{k,0};u^{k,0}) = (w^k;u^k)$. Therefore, PL-ADMM can be seen as a degenerate case of our IALM equipped with 0/1-BALM which only performs one step when solving subproblems in the setting of $t = 0$. However, the global convergence of PL-ADMM proposed in \citet{bolte2018nonconvex,boct2020proximal} requires that $A$ is nonsingular and the condition number of $A^\top A$ is less than 2. In the subsequent Theorems \ref{sub-fin-ter}, \ref{sub_con} and \ref{whole_con}, under a weaker assumption that $A$ is full row rank, our IALM equipped with 0/1-BALM could ensure the global convergence.
	\end{itemize}
\end{remark}



To proceed, we give the following basic assumptions.
\begin{assumption} \label{ass_fullrr}
	Matrix $A$ is full row rank with $\gamma : = \sqrt{\theta_{\min} (AA^\top)} > 0$.
\end{assumption}

\begin{assumption}\label{fbb}
	The function $f$ is bounded below.
\end{assumption}

Next let us show some convergence properties of 0/1-BALM. 

\begin{theorem}[Global Convergence of 0/1-BALM] \label{pam_converge} 
	Suppose that Assumptions \ref{bbs}, \ref{ass_fullrr} and \ref{fbb} hold, and parameters $\alpha$, $\{Q_{k,j}\}_{j \in \mathbb{N}}$ and $t$ satisfy
	\begin{align} \label{para_palm}
		\alpha = \frac{1}{\rho},~ \sup_{j \in \mathbb{N}} \| Q_{k,j} \| \leq q,~ t > q + l_g,
	\end{align}
	where $q$ is a positive constant. If $\{ (w^{k,j},u^{k,j}) \}_{j \in \mathbb{N}}$ is a sequence generated by 0/1-BALM, then the following assertions are true:
	
	(i) The objective function sequence satisfies the following descent property
	\begin{equation} \label{pam_suff_des}
		\mathcal{V}_{\rho,\mu} (w^{k,j},{u}^{k,j},z^k,v^{k}) - \mathcal{V}_{\rho,\mu} (w^{k,j+1},{u}^{k,j+1},z^k,v^{k}) \geq \frac{\zeta}{2} \| w^{k,j+1} - w^{k,j} \|^2 , \forall j \in \mathbb{N},
	\end{equation}
	where $\zeta: = t - q - l_g > 0$.
	
	(ii) We have 
	\begin{equation}\label{pam_j+1-j}
		\lim_{j \to \infty} \| w^{k,j+1} - w^{k,j} \| = 0,~ \lim_{j \to \infty} \| u^{k,j+1} - u^{k,j} \| = 0.
	\end{equation}
	
	(iii) The sequence $\{(w^{k,j},u^{k,j})\}_{j \in \mathbb{N}}$ is bounded.
	
	(iv) Each accumulation point of sequence $\{(w^{k,j},u^{k,j})\}_{j \in \mathbb{N}}$ is a ${\rm P}^k$-stationary point. Moreover, we have
	\begin{equation} \label{lim_nabh}
		\lim_{j \to \infty} \nabla_w h_k(w^{k,j},u^{k,j}) = 0.
	\end{equation}
	
	(v) If $f$ is a convex function, then the sequence $\{(w^{k,j},u^{k,j})\}_{j \in \mathbb{N}}$ converges to a ${\rm P}^k$-stationary point. 
	
\end{theorem}

\proof{Proof.}~(i) 
It follows from (\ref{u-step-sim1}) and the definition of proximal operator, we have
\begin{equation} \notag \label{u_stepeq}
	\frac{1}{2} \| Aw^{k,j} + b - u^{k,j+1} + z^k/ \rho \|^2 + (\lambda/\rho) \| u^{k,j+1}_+ \|_0 \leq \frac{1}{2} \| Aw^{k,j} + b - u^{k,j} + z^k/ \rho \|^2 + (\lambda/\rho) \| u^{k,j}_+ \|_0.
\end{equation}
Multiplying $\rho$ and adding $f(w^{k,j}) + (\mu/2) \| w^{k,j} - v^{k} \|^2 - (1/2\rho) \| z^k \|^2$ on both side of above inequality, we can obtain
\begin{align}
	&\mathcal{V}_{\rho,\mu} (w^{k,j},{u}^{k,j},z^k,v^{k}) - \mathcal{V}_{\rho,\mu} (w^{k,j},{u}^{k,j+1},z^k,v^{k}) \geq 0. \label{suff_des_u}
\end{align}
From the definition of $\mathcal{V}_{\rho,\mu}$ and $\varphi_{k,j}$, the following equations hold for any $j \in \mathbb{N}$
\begin{align} &\label{v_htilde}
	\mathcal{V}_{\rho,\mu} (w,{u}^{k,j+1},z^k,v^{k}) \notag \\
	=& \varphi_{k,j}(w) - \mathcal{D}_{{\psi }_{k,j}}(w, w^{k,j}) + \lambda \| u^{k,j+1}_+ \|_0 \notag \\
	\mathop{=}\limits^{(\ref{breg_sim})}& \varphi_{k,j}(w) + g(w) - \langle \nabla g(w^{k,j}), w - w^{k,j} \rangle - \frac{1}{2} \| w - w^{k,j} \|^2_{P_{k,j}} + \lambda \| u^{k,j+1}_+ \|_0 - g(w^{k,j}).
\end{align}
This implies that
\begin{align}
	&\mathcal{V}_{\rho,\mu} (w^{k,j},{u}^{k,j+1},z^k,v^k) - \mathcal{V}_{\rho,\mu} (w^{k,j+1},{u}^{k,j+1},z^k,v^k) \notag \\ 
	= & \frac{1}{2} \| w^{k,j+1} - w^{k,j} \|^2_{P_{k,j}} - ( g(w^{k,j+1}) - g(w^{k,j}) - \langle \nabla g(w^{k,j}), w^{k,j+1} - w^{k,j} \rangle ) \notag \\ 
	&+ \varphi_{k,j}(w^{k,j}) - \varphi_{k,j}(w^{k,j+1}) \label{w-descent} \\
	\mathop{\geq}\limits^{(\ref{w-step})} & \frac{1}{2} \| w^{k,j+1} - w^{k,j} \|^2_{P_{k,j}} - ( g(w^{k,j+1}) - g(w^{k,j}) - \langle \nabla g(w^{k,j}), w^{k,j+1} - w^{k,j} \rangle )  \notag \\
	\geq & \frac{t - q -  l_g }{2} \| w^{k,j+1} - w^{k,j} \|^2 = \frac{\zeta}{2} \| w^{k,j+1} - w^{k,j} \|^2, \label{w-step1}	
\end{align}
where the second inequality follows from Lemma \ref{dl} and (\ref{para_palm}). Then adding (\ref{w-step1}) and (\ref{suff_des_u}) yields that
\begin{equation} \notag
	\mathcal{V}_{\rho,\mu} (w^{k,j},{u}^{k,j},z^k,v^{k}) - \mathcal{V}_{\rho,\mu} (w^{k,j+1},{u}^{k,j+1},z^k,v^{k}) \geq \frac{\zeta}{2}\| w^{k,j+1} - w^{k,j} \|^2 .
\end{equation}

(ii) By the definition of Lyapunov function, we have
\begin{align}
	\mathcal{V}_{\rho,\mu} (w,{u},z^k,v^{k}) &= f(w) + \lambda \| u_+ \|_0 + \langle z^k, Aw + b - u \rangle + \frac{\rho}{2} \| Aw + b - u \|^2 + \frac{\mu}{2} \| w - v^{k} \|^2 \notag \\
	& = f(w) + \lambda \| u_+ \|_0  + \frac{\rho}{2} \| Aw + b - u + \frac{1}{\rho} z^k \|^2 - \frac{1}{2\rho} \| z^k \|^2 + \frac{\mu}{2} \| w - v^{k} \|^2. \notag
\end{align} 
Notice that $f(w)$ is bounded below from Assumption \ref{fbb}. Then  $\mathcal{V}_{\rho,\mu} (w,{u},z^k,v^{k})$ is bounded below when $z^k$ and $v^k$ are fixed. Moreover, taking (\ref{pam_suff_des}) into account, $\{ \mathcal{V}_{\rho,\mu} (w^{k,j},u^{k,j},z^k,v^{k}) \}_{j \in \mathbb{N}}$ is a decreasing and bounded below sequence, which means that 
\begin{equation} \notag
	\lim_{j \to \infty} \mathcal{V}_{\rho,\mu} (w^{k,j},u^{k,j},z^k,v^{k}) = \widehat{\mathcal{V}},
\end{equation}
where $\widehat{\mathcal{V}}$ is a finite real number.
Then taking the limit in (\ref{pam_suff_des}) as $j \to \infty$, we can obtain
\begin{equation}  \label{succ_w}
	\lim_{j \to \infty} \| w^{k,j+1} - w^{k,j} \| = 0.
\end{equation}
We can deduce from (\ref{opt_phi}) that
\begin{align}
	\rho A^\top( u^{k,j+1} - u^{k,j} ) =& \nabla f( w^{k,j+1} ) - \nabla f( w^{k,j} ) + ( \mu I + \rho A^\top A + P_{k,j} ) ( w^{k,j+1} - w^{k,j} ) \notag \\
	& - P_{k,j-1}( w^{k,j} - w^{k,j-1} ) - (\nabla g(w^{k,j+1}) - \nabla g(w^{k,j}) ) +  (\nabla g(w^{k,j}) - \nabla g(w^{k,j-1}) ). \notag
\end{align}
Then we have the following estimation
\begin{align}
	&\rho \gamma \| u^{k,j+1} - u^{k,j} \| \notag  
	\leq \rho \| A^\top( u^{k,j+1} - u^{k,j} ) \| \notag \\
	\leq & \| \nabla f( w^{k,j+1} ) - \nabla f( w^{k,j} ) \| + (\mu+ \rho\| A \|^2 + \| P_{k,j} \|) \| w^{k,j+1} - w^{k,j} \| + \| P_{k,j-1} \| \| w^{k,j} - w^{k,j-1} \| \notag \\
	&+ \| \nabla g(w^{k,j+1}) - \nabla g(w^{k,j}) \| + \| \nabla g(w^{k,j}) - \nabla g(w^{k,j-1}) \| \notag \\
	\leq & (\mu+ \rho\| A \|^2 + t + l_g+ l_f) \| w^{k,j+1} - w^{k,j} \| + (t + l_g ) \| w^{k,j} - w^{k,j-1} \|, \label{u_upper}
\end{align}
where the first inequality follows from Assumption \ref{ass_fullrr}, the last inequality follows from $\| P_{k,j} \| \leq t$, the $l_f$-smoothness of $f$ and the $l_g$-smoothness of $g$.
Since (\ref{succ_w}) holds, taking limit as $j \to \infty$ in (\ref{u_upper}), we can obtain
\begin{equation}  \notag
	\lim_{j \to \infty} \| u^{k,j+1} - u^{k,j} \| = 0.
\end{equation}
We arrive at the desire assertion (ii).

(iii) By the decreasing property of $\{ \mathcal{V}_{\rho,\mu} (w^{k,j},u^{k,j},z^k,v^{k}) \}_{j \in \mathbb{N}}$, we have
\begin{align} \notag \label{pam_upbound}
	&\mathcal{V}_{\rho,\mu} (w^{k,0},u^{k,0},z^k,v^{k}) \\ 
	\geq& \mathcal{V}_{\rho,\mu} (w^{k,j},u^{k,j},z^k,v^{k}) \notag \\
	= & f(w^{k,j}) + \lambda \| u^{k,j}_+ \|_0 +  \frac{\rho}{2} \| Aw^{k,j} + b - u^{k,j} + \frac{1}{\rho} z^k \|^2 - \frac{1}{2\rho} \| z^k \|^2 + \frac{\mu}{2} \| w^{k,j} - v^{k} \|^2 \notag \\
	\geq & \inf_{ w \in \mathbb{R}^p } f(w) + \lambda \| u^{k,j}_+ \|_0 +  \frac{\rho}{2} \| Aw^{k,j} + b - u^{k,j} + \frac{1}{\rho} z^k \|^2 - \frac{1}{2\rho} \| z^k \|^2 + \frac{\mu}{2} \| w^{k,j} - v^{k} \|^2, \notag 
\end{align}
where the last inequality follows from the lower boundedness of $f$. We can further obtain 
\begin{align}
	&\mathcal{V}_{\rho,\mu} (w^{k,0},u^{k,0},z^k,v^{k}) - \inf_{ w \in \mathbb{R}^p } f(w) + \frac{1}{2\rho} \| z^k \|^2 \notag \\
	\geq& \lambda \| u^{k,j}_+ \|_0 + \frac{\rho}{2} \| Aw^{k,j} + b - u^{k,j} - \frac{1}{\rho} z^k \|^2 + \frac{\mu}{2} \| w^{k,j} - v^{k} \|^2. \notag 
\end{align}
This implies that $\{ w^{k,j} - v^{k} \}_{j \in \mathbb{N}}$ and $\{ Aw^{k,j} + b - u^{k,j} - z^k / \rho \}_{j \in \mathbb{N}}$ is bounded. Then it follows from $ \| w^{k,j} \| \leq \| v^{k} \| + \| w^{k,j} - v^{k} \|$ and $ \| u^{k,j} \| \leq \| Aw^{k,j} + b - u^{k,j} - z^k / \rho \| + \| A \| \| w^{k,j} \| + \| b \| +  \| z^k \| / \rho $ that $\{ w^{k,j} \}_{j \in \mathbb{N}}$ and $\{ u^{k,j} \}_{j \in \mathbb{N}}$ is bounded. We complete the proof of assertion (iii).

(iv) The assertion (iii) implies that $\{ (w^{k,j};u^{k,j}) \}_{j \in \mathbb{N}}$ has an accumulation point. Suppose $(\widehat{w},\widehat{u})$ is an accumulation point of the sequence $\{ (w^{k,j},u^{k,j}) \}_{j \in \mathbb{N}}$, which means that there exists an infinite set $J \subseteq \mathbb{N}$ such that 
\begin{equation} \notag \label{wu-sublimit}
	\lim_{\substack{j \to \infty \\ j \in J}} w^{k,j} = \widehat{w},~\lim_{\substack{j \to \infty \\ j \in J}} u^{k,j} = \widehat{u}.
\end{equation}
Then from (\ref{pam_j+1-j}), we have
\begin{subequations}  \label{prox_lim} \notag
	\begin{align}
		&\lim_{\substack{j \to \infty \\ j \in J}} u^{k,j+1} = \lim_{\substack{j \to \infty \\ j \in J}} (u^{k,j+1} - u^{k,j} +u^{k,j}) = \widehat{u}, \\
		& \lim_{\substack{j \to \infty \\ j \in J}} (u^{k,j} - \alpha \nabla_u h_k(w^{k,j},u^{k,j})) = \widehat{u} - \alpha \nabla_u h_k(\widehat{w},\widehat{u}). 
	\end{align}
\end{subequations}
Furthermore, taking (\ref{u-step}) into account, the proximal behavior (see e.g. \citet[Theorem 1.25]{RockWets98}) implies that
\begin{equation} \label{P-sta-u}
	\widehat{u} \in {\rm Prox}_{\alpha \lambda \| (\cdot)_+ \|_0} ( \widehat{u} - \alpha \nabla_u h_k(\widehat{w},\widehat{u}) ).
\end{equation}
From (\ref{opt_phi}), we can obtain
\begin{align} 
	&\nabla f(w^{k,j}) + \mu (w^{k,j} - v^{k}) + A^\top (z^k + \rho( Aw^{k,j} + b - u^{k,j} )) + \nabla f(w^{k,j+1}) - \nabla f(w^{k,j}) \notag \\
	& + (\mu I + \rho A^\top A + P_{k,j}) ( w^{k,j+1} - w^{k,j} ) - ( \nabla g(w^{k,j+1}) - \nabla g(w^{k,j}) ) - \rho A^\top (u^{k,j+1} - u^{k,j}) = 0.\label{opt-w-step} \notag
\end{align}
Taking limit as $j \to \infty$ for $j \in J$ in above equality, it follows from (\ref{pam_j+1-j}), and the Lipschitz continuity of $f$ and $g$ that
\begin{equation} \notag \label{P-sta-w}
	\nabla f(\widehat{w}) + \mu (\widehat{w} - v^{k}) + A^\top (z^k + \rho( A\widehat{w} + b - \widehat{u} )) = 0.
\end{equation}
Taking $\widehat{z} = - \nabla_u h_k( \widehat{w},\widehat{u} ) = z^k + \rho (A\widehat{w} + b - \widehat{u})$, then the above equality and (\ref{P-sta-u}) can be rewritten as
\begin{align}
	& \nabla f(\widehat{w}) + \mu (\widehat{w} - v^{k}) + A^\top \widehat{z} = 0, \label{pks_con1}\\
	& \widehat{u} \in {\rm Prox_{\alpha\lambda \|(\cdot)_+\|_0}} ( \widehat{u} + \alpha \widehat{z} ), \notag
\end{align}
which means that $(\widehat{w};\widehat{u})$ is a ${\rm P}^k$-stationary point.

Denote $\Theta$ as the set including all the accumulation points of sequence $\{(w^{k,j};u^{k,j})\}_{j \in \mathbb{N}}$. It follows from \citet[Lemma 5 (ii)]{bolte2014proximal} that 
\begin{equation} \label{distpomeg}
	\lim_{j \to \infty} {\rm dist} ( (w^{k,j};u^{k,j}); \Theta ) = 0.
\end{equation}
Taking $\forall (\widehat{w};\widehat{u}) \in \Theta$, then $(\widehat{w};\widehat{u})$ is a ${\rm P}^k$-stationary point and we can calculate
\begin{align}
	\nabla_w h_k( \widehat{w},\widehat{u} ) & = \nabla f(\widehat{w}) + \mu  (\widehat{w} - v^{k} ) + \rho A^\top ( A\widehat{w} + b -\widehat{u} + \frac{1}{\rho} z^k ) \notag \\
	& = \nabla f(\widehat{w}) + \mu (\widehat{w} - v^{k}) + A^\top \widehat{z} \mathop{=}\limits^{(\ref{pks_con1})} 0. \label{nabwh0}
\end{align}
We have the following estimation
\begin{align}
	\| \nabla_w h_k( w^{k,j}, u^{k,j} ) \| \mathop{=}\limits^{(\ref{nabwh0})}& \| \nabla_w h_k( w^{k,j}, u^{k,j} ) - \nabla_w h_k( \widehat{w},\widehat{u} ) \| \leq  l_h \| ( w^{k,j} - \widehat{w}; u^{k,j} - \widehat{u} ) \|, \label{ub_gh}
\end{align}
where the  last equality follows from Property \ref{prop_hstar} (ii). Taking limit inferior for all $(\widehat{w};\widehat{u}) \in \Theta$ in right side (\ref{ub_gh}), we can obtain
\begin{equation} \notag \label{nabh_ub1}
	\| \nabla_w h_k( w^{k,j}, u^{k,j} ) \| \leq l_h {\rm dist} ( ( w^{k,j}; u^{k,j}) ; \Theta ) .
\end{equation}
Note that (\ref{distpomeg}) holds. Taking limit as $j \to \infty$ in above inequality, we can obtain the desired inclusion.

(v) Let us first prove that there are only finite ${\rm P}^k$-stationary points for (\ref{eq2.1}). Suppose that $(\widehat{w};\widehat{u})$ is a ${\rm P}^k$-stationary point with constant $\alpha > 0$ and $\widehat{z}$ is the corresponding ${\rm P}^k$-stationary multiplier, from (\ref{P^k-stat}), $(\widehat{w};\widehat{u};\widehat{z})$ satisfies
\begin{equation}   \label{Pk_stat_eq}
	\left\{ \begin{aligned}
		& \nabla_w h_k( \widehat{w}, \widehat{u}) = 0,\\
		& \widehat{u} \in {\rm Prox_{\alpha\lambda \|(\cdot)_+\|_0}} ( \widehat{u} + \alpha \widehat{z} ), \\
		& \widehat{z} = -\nabla_u h_k( \widehat{w}, \widehat{u}).
	\end{aligned}  \right.
\end{equation}
Moreover, given the vectors $\widehat{u}$ and $\widehat{z}$ satisfying the second formula of (\ref{Pk_stat_eq}), taking the same procedure as Property \ref{rela_cr_p}, we know that (\ref{uz_prox}) holds. Denote $\widehat{S}_- : = \{ i \in [n]: \widehat{u}_i \leq 0 \}$. Let us consider the following convex programming 
\begin{equation}
	\begin{aligned} \label{sub-conprog}
		&\min\limits_{w \in \mathbb{R}^p, u \in \mathbb{R}^n} && h_k(w,u) \\
		&~~~~~~s.t. && u_{\widehat{S}_-} \leq 0.	
	\end{aligned}
\end{equation}
Since $f$ is convex, Property \ref{prop_hstar} (iii) implies that $h_k(w,u)$ is a strongly convex function. Moreover, the constraints are linear inequalities. Therefore, (\ref{sub-conprog}) has the unique global minimizer. The KKT condition of (\ref{sub-conprog}) can be represented as
\begin{equation} \label{sub_kkt}
	\left\{  \begin{aligned}
		& \nabla_w h_k(w,u) = 0,~[\nabla_u h_k(w,u)]_{\widehat{S}_-^c} = 0, \\
		& z_{\widehat{S}_-} \geq 0 ,~ u_{\widehat{S}_-} \leq 0,~\langle z_{\widehat{S}_-},u_{\widehat{S}_-}  \rangle = 0,\\
		& z_{\widehat{S}_-} = -[\nabla_u h_k(w,u)]_{\widehat{S}_-},
	\end{aligned}  \right.
\end{equation}
where $z_{\widehat{S}_-} \in \mathbb{R}^{| \widehat{S}_- |}$ is the Lagrange multiplier. From (\ref{uz_prox}) and (\ref{Pk_stat_eq}), we know that $(\widehat{w};\widehat{u})$ satisfies (\ref{sub_kkt}). Then by the property of convex programming, we can conclude that $(\widehat{w};\widehat{u})$ is the unique global minimizer of (\ref{sub-conprog}). Notice that there are only a finite number of $\widehat{S}_-$ (no more than $2^n$). Therefore, the number of ${\rm P}^k$-stationary points for (\ref{eq2.1}) is finite. Because we have proved that each accumulation point of $\{ (w^{k,j},u^{k,j}) \}_{j \in \mathbb{N}}$ is a ${\rm P}^k$-stationary point, the number of accumulation points is finite and hence each accumulation point of $\{ (w^{k,j},u^{k,j}) \}_{j \in \mathbb{N}}$ is isolated. Considering that (\ref{pam_j+1-j}) holds, it follows from \citet[Proposition 7]{kanzow1999qp} that the whole sequence $\{ (w^j,u^j) \}_{j \in \mathbb{N}}$ must converge to a ${\rm P}^k$-stationary point $(\widehat{w};\widehat{u})$. 
\hfill \Halmos
\endproof

Next we will show that the 0/1-BALM satisfies the stopping criteria (\ref{stopping1}) within finite iterations. Thus our IALM is well defined.

\begin{theorem} \label{sub-fin-ter}
	Suppose that Assumptions \ref{bbs}, \ref{ass_fullrr} and \ref{fbb} hold, and parameters $ \alpha, \{Q_{k,j}\}_{j \in \mathbb{N}}$ and $t$ satisfy (\ref{para_palm}). If $\{ (w^{k,j};u^{k,j}) \}_{j \in \mathbb{N}}$ is a sequence generated by 0/1-BALM,
	then there exists $\overline{j} \in \mathbb{N} $ such that $(w^{k+1};u^{k+1}) = (w^{k,\overline{j} + 1};u^{k,\overline{j} + 1})$ satisfies the stopping criteria (\ref{stopping1}).
\end{theorem}

\proof{Proof.}
Note that (\ref{pam_suff_des}) holds true. Therefore, for any $j \in \mathbb{N}$, taking $(w^{k+1};u^{k+1}) = (w^{k,j+1};u^{k,j+1})$ leads to the stopping criteria (\ref{des_auglag}). By (\ref{u-step}), we have the following estimation for all $j \in \mathbb{N}$
\begin{equation} \label{dist_upbound}
	{\rm dist}( u^{k,j}, {\rm Prox_{\alpha\lambda \|(\cdot)_+\|_0}} ( u^{k,j} - \alpha \nabla_u h_k ( w^{k,j},u^{k,j} ) ) ) \leq \| u^{k,j} - u^{k,j+1} \|. \notag
\end{equation}
Note that (\ref{pam_j+1-j}) holds true. Taking limit as $j \to \infty$ in above inequality leads to
\begin{equation} \notag \label{dist_lim}
	\lim_{j \to \infty} {\rm dist}( u^{k,j}, {\rm Prox_{\alpha\lambda \|(\cdot)_+\|_0}} ( u^{k,j} - \alpha \nabla_u h_k ( w^{k,j},u^{k,j} ) ) ) = 0.
\end{equation}	
This and (\ref{lim_nabh}) yield that for any $\epsilon_k > 0$, there exists $\overline{j} \in \mathbb{N}$ such that the following inequality holds
\begin{align} \notag
	\max \{ \| \nabla_w h_k ( w^{k,\overline{j} + 1}, u^{k,\overline{j}+1} ) \|, {\rm dist}( u^{k,\overline{j} + 1}, {\rm Prox_{\alpha\lambda \|(\cdot)_+\|_0}} ( u^{k,\overline{j} + 1} - \alpha \nabla_u h_k ( w^{k,\overline{j} + 1},u^{k,\overline{j} + 1} ) ) ) \} \leq \epsilon_k.
\end{align}
Overall, we obtain the desired conclusion. \hfill\Halmos
\endproof 

At the end of this section, let us discuss the choice of parameter $t$ for the two special cases in Remark \ref{re-palm}.

\begin{remark} For some special cases, we could give a more intuitive and larger range of $t$ compared with the setting (\ref{para_palm}) in Theorems \ref{pam_converge} and \ref{sub-fin-ter}. Actually, it suffices to choose parameter $t$ satisfying the sufficient descent result in Theorem \ref{pam_converge} (i).
	
	For the case of $Q_{k,j} = \rho A^\top A$ and $g(w) = f(w)$, we can choose $t > (l_f+\rho \| A \|^2 - \mu)/2$.  Note that in this setting, $\varphi_{k,j}$ is a quadratic function. Then we have the following estimation by the Taylor expansion,
	\begin{align} \notag
		\varphi_{k,j} ( w^{k,j} ) - \varphi_{k,j} ( w^{k,j+1} ) & = \langle \nabla \varphi_{k,j} ( w^{k,j+1} ), w^{k,j} - w^{k,j+1} \rangle + \frac{t + \mu}{2} \| w^{k,j+1} - w^{k,j} \|^2 \notag \\
		& \mathop{=}\limits^{(\ref{w-step})} \frac{t+\mu}{2} \| w^{k,j+1} - w^{k,j} \|^2.  \notag
	\end{align}
	This and (\ref{w-descent}) imply that 
	\begin{align}
		&\mathcal{V}_{\rho,\mu} (w^{k,j},{u}^{k,j+1},z^k,v^{k}) - \mathcal{V}_{\rho,\mu} (w^{k,j+1},{u}^{k,j+1},z^k,v^{k}) \notag \\ 
		\mathop{=}\limits^{} & \frac{1}{2} \| w^{k,j+1} - w^{k,j} \|^2_{tI - \rho A^\top A} - ( f(w^{k,j+1}) - f(w^{k,j}) - \langle \nabla f(w^{k,j}), w^{k,j+1} - w^{k,j} \rangle ) +  \frac{t+\mu}{2} \| w^{k,j+1} - w^{k,j} \|^2 \notag  \\
		\geq & \frac{2t + \mu -  l_f - \rho \| A \|^2 }{2} \| w^{k,j+1} - w^{k,j} \|^2, \notag
	\end{align}
	where $(2t + \mu -  l_f - \rho \| A \|^2) /2 > 0$ and hence Theorem \ref{pam_converge} (i) holds.
	
	For the case of $Q_{k,j} = \rho A^\top_{T_{k,j}^c} A_{T_{k,j}^c}$, $g(w) = 0$ and $f(w) = \| w \|^2/2$, we can choose $t > (\rho \| A \|^2 - \mu - 1)/2$. In this setting, $\phi_{k,j}$ is a quadratic function and its Taylor expansion can be represented as 
	\begin{align}
		\varphi_{k,j} ( w^{k,j} ) - \varphi_{k,j} ( w^{k,j+1} ) & = \langle \nabla \varphi_{k,j} ( w^{k,j+1} ), w^{k,j} - w^{k,j+1} \rangle + \frac{1}{2} \| w^{k,j+1} - w^{k,j} \|^2_{(t+\mu + 1)I + \rho A^\top_{T_{k,j}}A_{T_{k,j}}} \notag \\
		& \mathop{=}\limits^{(\ref{w-step})} \frac{1}{2} \| w^{k,j+1} - w^{k,j} \|^2_{(t+\mu + 1)I + \rho A^\top_{T_{k,j}}A_{T_{k,j}}}.  \notag
	\end{align}
	This together with (\ref{w-descent}) imply that 
	\begin{align}
		&\mathcal{V}_{\rho,\mu} (w^{k,j},{u}^{k,j+1},z^k,v^{k}) - \mathcal{V}_{\rho,\mu} (w^{k,j+1},{u}^{k,j+1},z^k,v^{k}) \notag \\ 
		\mathop{=}\limits^{} & \frac{1}{2} \| w^{k,j+1} - w^{k,j} \|^2_{tI-\rho A^\top_{T_{k,j}^c} A_{T_{k,j}^c}} +  \frac{1}{2} \| w^{k,j+1} - w^{k,j} \|^2_{(t+\mu + 1)I + \rho A^\top_{T_{k,j}}A_{T_{k,j}}} \notag \\
		\geq & \frac{2t+\mu+1 - \rho \|A\|^2}{2} \| w^{k,j+1} - w^{k,j} \|^2, \notag
	\end{align}
	where $(2t+\mu+1 - \rho \|A\|^2)/2 > 0$ and hence Theorem \ref{pam_converge} (i) holds.
\end{remark}

\section{Convergence Analysis of IALM.} \label{convergence_analysis}
In this section, we will establish the global convergence of IALM. Before this, let us present the setting of parameters and give a basic assumption.

Given parameter $\mu > 0$, the setting of parameters $\alpha$, $\rho$, $\beta$ and $\eta$ are as follows:
\begin{align} 
	& \rho > \max\left\{ \frac{16(c_1^2+c_2^2)}{\mu },  \frac{6l_f}{\gamma^2} \right\},~\alpha = \frac{1}{\rho},~ \beta := \frac{8c_2^2}{\rho},~ \eta >  \frac{4}{\rho \gamma^2}, \label{para_mu}
\end{align}
where $c_1:=(\mu+l_f)/\gamma, c_2: = \mu/\gamma$. Moreover, the sequence $\{\epsilon_k\}_{k \in \mathbb{N}}$ satisfies 
\begin{equation} \label{epsilon_val}
	\epsilon_{k+1} / \epsilon_k \leq \sqrt{ (\rho\gamma^2\eta - 4) / (\rho\gamma^2\eta + 4)}.
\end{equation}
For convenience, we denote $\mathcal{V}_k := \mathcal{V}_{\rho,\beta}( w^k,u^k,z^k,v^{k-1} )$ for all $k \in \mathbb{N}^+$. Note that we adopt the penalty parameter $\beta$ here instead of the $\mu$ in Algorithm \ref{IALM}. Moreover, we take the Lyapunov variable as $v^{k-1}$. Let us give the sufficient decrease property of sequence $\{ \mathcal{V}_k + \eta \epsilon_{k-1}^2 \}_{k\in \mathbb{N}^+}$.  
\begin{lemma}[Descent Property] \label{lemma3.1}Suppose that Assumptions \ref{bbs}, \ref{ass_fullrr} and \ref{fbb} hold, and parameters are chosen as (\ref{para_mu}) and (\ref{epsilon_val}). Let $\{ (w^k,u^k,z^k,v^k) \}_{k\in \mathbb{N}}$ be a sequence generated by algorithm IALM satisfying stopping criteria (\ref{stopping1}). Then we have
	\begin{equation}  \label{lyap_des}
		(\mathcal{V}_k + \eta \epsilon_{k-1}^2) - (\mathcal{V}_{k+1} + \eta \epsilon_{k}^2) \geq \tau \| w^{k+1} - w^k \|^2,
	\end{equation}
	where $\tau:= \mu/4$.
\end{lemma}

\proof{Proof.}
From (\ref{eq2.2}), (\ref{nabla_h}) and $v^k = w^k$, we have
\begin{align} 
	&\nabla_w h_k (w^{k+1},u^{k+1}) - \nabla_w h_{k-1} (w^{k},u^{k}) \notag \\
	=& \nabla f(w^{k+1}) - \nabla f(w^k) + \mu (w^{k+1} - w^k) - \mu (w^k - w^{k-1}) + A^\top (z^{k+1} - z^k). \label{eq3.4}  \notag
\end{align}
From above equality, (\ref{crit1}) and Assumption \ref{ass_fullrr}, we can deduce
\begin{align} 
	\gamma \|z^{k+1} - z^k \| \leq & \| A^\top ( z^{k+1} - z^k )  \| \notag \\
	\leq & (\mu+l_f) \| w^{k+1} - w^k \| + \mu \| w^k - w^{k-1} \|  + \|\nabla_w h_{k-1}(w^k,u^k)\| + 	\|\nabla_w h_{k}(w^{k+1},u^{k+1})\| \notag\\
	\leq & (\mu+l_f) \| w^{k+1} - w^k \| + \mu \| w^k - w^{k-1} \| + 	\epsilon_{k-1} + 	\epsilon_{k}, \notag
\end{align}
which directly leads to
\begin{equation} \label{eq3.5}
	\| z^{k+1} - z^k \| \leq c_1 \| w^{k+1} - w^k \| + c_2 \| w^k - w^{k-1} \| + \epsilon_{k-1}/\gamma + \epsilon_{k}/\gamma.
\end{equation}
Furthermore, using the fact $(\sum_{i=1}^4 t_i)^2 \leq 4\sum_{i=1}^4 t_i^2 $, we can estimate from (\ref{eq3.5}) that
\begin{equation} \label{eq3.6}
	\| z^{k+1} - z^k \|^2 \leq 4c_1^2 \| w^{k+1} - w^k \|^2 + 4c_2^2 \| w^k - w^{k-1} \|^2 + \frac{4}{\gamma^2} \epsilon_{k-1}^2 + \frac{4}{\gamma^2} \epsilon_{k}^2 .
\end{equation}
From (\ref{eq4}), $v^{k-1} = w^{k-1}$ and $v^k = w^k$, (\ref{des_auglag}) yields that
\begin{equation} \label{auglag_des}
	\mathcal{L}_\rho ( w^k,u^k,z^k ) - \mathcal{L}_\rho ( w^{k+1},u^{k+1},z^k ) \geq \frac{\mu}{2} \| w^{k+1} - w^{k} \|^2.
\end{equation}
Then we have the following estimation
\begin{align}
	&\mathcal{L}_\rho( w^k,u^k,z^k ) - \mathcal{L}_\rho( w^{k+1},u^{k+1},z^{k+1} ) \nonumber \\
	=& \mathcal{L}_\rho( w^k,u^k,z^k )  - \mathcal{L}_\rho( w^{k+1},u^{k+1},z^k ) + \mathcal{L}_\rho( w^{k+1},u^{k+1},z^k ) - \mathcal{L}_\rho( w^{k+1},u^{k+1},z^{k+1} ) \nonumber \\
	\mathop{\geq}\limits^{(\ref{auglag_des})}& \frac{\mu}{2} \| w^{k+1} - w^k \|^2 + \langle z^k - z^{k+1}, Aw^{k+1} + b - u^{k+1} \rangle \notag \\
	\mathop{\geq}\limits^{(\ref{eq2.2})}&\frac{\mu}{2} \| w^{k+1} - w^k \|^2 - \frac{1}{\rho} 	\| z^{k+1} - z^k \|^2 \notag \\
	\mathop{\geq}\limits^{(\ref{eq3.6})}& ( \frac{\mu}{2} - \frac{4c_1^2}{\rho}) \| w^{k+1} - w^k \|^2 - \frac{4c_2^2}{\rho}  \| w^k - w^{k-1} \|^2 - \frac{4}{\rho\gamma^2} \epsilon_{k-1}^2 - \frac{4}{\rho\gamma^2} \epsilon_{k}^2. \label{auglag_diff}
\end{align}
From (\ref{eq4}) and $v^k = w^k$, we can obtain
\begin{equation} \notag
	\begin{aligned} 
		& \mathcal{V}_k + \eta \epsilon_{k-1}^2 - (\mathcal{V}_{k+1} + \eta \epsilon_k^2) \\
		=& \mathcal{L}_\rho( w^k,u^k,z^k ) - \mathcal{L}_\rho( w^{k+1},u^{k+1},z^{k+1} ) + \frac{\beta}{2}\big( \| w^{k} - w^{k-1} \|^2 - \| w^{k+1} - w^k \|^2 \big) + \eta \epsilon_{k-1}^2 - \eta \epsilon_k^2 \notag\\
		\mathop{\geq}\limits^{(\ref{auglag_diff})}& \big( \frac{\mu}{2} - \frac{4c_1^2}{\rho} - \frac{\beta}{2} \big) \| w^{k+1} - w^k \|^2 + \big( \frac{\beta}{2} - \frac{4c_2^2}{\rho} \big) \| w^k - w^{k-1} \|^2 + [(\eta - \frac{4}{\rho\gamma^2})  \epsilon_{k-1}^2 - (\eta + \frac{4}{\rho\gamma^2})  \epsilon_{k}^2]  \\
		\mathop{\geq}\limits^{(\ref{para_mu},\ref{epsilon_val})}& \tau  \| w^{k+1} - w^k \|^2,
	\end{aligned}
\end{equation}
which means that the descent property of sequence $\{ \mathcal{V}_k + \eta \epsilon_{k-1}^2 \}_{k\in \mathbb{N}^+}$  holds. \hfill\Halmos
\endproof

To ensure the boundedness of the sequence, we need the following assumption.

\begin{assumption} \label{assum3.1}
	The function $f$ is coercive.
\end{assumption}

Actually, for the smooth function $f$, Assumption \ref{assum3.1} is stronger than Assumption \ref{fbb}. The next lemma shows that the sequence generated by IALM is bounded.

\begin{lemma}[Sequence Boundedness] \label{seq_bounded}Suppose that Assumptions \ref{bbs}, \ref{ass_fullrr} and \ref{assum3.1} hold, and parameters are selected satisfied (\ref{para_mu}) and (\ref{epsilon_val}). Let $\{ (w^k;u^k;z^k;v^k) \}_{k\in \mathbb{N}}$ be a sequence generated by IALM satisfying stopping criteria (\ref{stopping1}). Then $\{ (w^k;u^k;z^k;v^k) \}_{k\in \mathbb{N}}$ is bounded.
\end{lemma}

\proof{Proof.} 
From (\ref{crit1}), (\ref{eq2.2}), $v^k = w^k$ and Assumption \ref{ass_fullrr}, we can estimate an upper bound for $\{z^k\}_{k\in \mathbb{N}}$
\begin{equation} \notag
	\gamma \left\|z^{k+1}\right\| \leq  \|A^\top z^{k+1}\| \leq \| \nabla f(w^{k+1}) \| + \mu \| w^{k+1} - w^{k} \| + \epsilon_{k} ,
\end{equation}
which means that
\begin{equation} \label{z_bound1}
	\left\|z^{k+1}\right\| \leq \frac{1}{\gamma} \| \nabla f(w^{k+1}) \| + \frac{\mu}{\gamma} \| w^{k+1} - w^{k} \| + \frac{\epsilon_{k}}{\gamma}.
\end{equation}
Then using the fact $(\sum_{i=1}^3 t_i)^2 \leq 3\sum_{i=1}^3 t_i^2 $ for real number $t_i, i=1,2,3$, we can obtain
\begin{equation} \label{z_bound}
	\|z^{k+1}\|^2 \leq \frac{3}{\gamma^2} \| \nabla f(w^{k+1}) \|^2 + \frac{3\mu^2}{\gamma^2} \| w^{k+1} - w^{k} \|^2 + \frac{3\epsilon_{k}^2}{\gamma^2}.
\end{equation} 
Since $f$ is $l_f$-smooth, we can give the following estimation by Lemma \ref{dl}
\begin{align}
	&f( w^{k+1} - \frac{1}{l_f} \nabla f( w^{k+1} ) ) \notag \\ \leq & f(w^{k+1}) + \langle \nabla f(w^{k+1}) , w^{k+1} - \frac{1}{l_f} \nabla f(w^{k+1}) - w^{k+1} \rangle  +  \frac{l_f}{2} \| w^{k+1} - \frac{1}{l_f} \nabla f(w^{k+1}) - w^{k+1} \|^2 \notag \\
	= & f(w^{k+1}) - \frac{1}{2l_f} \| \nabla f(w^{k+1}) \|^2. \label{f-nabla_f}
\end{align}
Noting that $f$ is smooth and coercive, we know that $f$ is bounded below. 
Considering the conclusion of Lemma \ref{lemma3.1}, we have the following inequalities for any $k \in \mathbb{N}^+$, 
\begin{align}
	&~~~\mathcal{V}_1 + \eta \epsilon_0^2 \notag \\
	&\geq (\mathcal{V}_{k+1} + \eta \epsilon_{k}^2 ) + (\mathcal{V}_k + \eta \epsilon_{k-1}^2 ) - (\mathcal{V}_{k+1} + \eta \epsilon_{k}^2 ) \geq \mathcal{V}_{k+1} + \eta \epsilon_{k}^2 + \tau \| w^{k+1} - w^k \|^2  \notag \\
	&= f(w^{k+1}) + \lambda \|u^{k+1}_+\|_0 + \frac{\rho}{2} \| Aw^{k+1} + b - u^{k+1} + \frac{1}{\rho}z^{k+1} \|^2 \notag \\
	&~~~ + (\tau + \frac{\beta}{2})\| w^{k+1} - w^k \|^2+ \eta \epsilon_{k}^2 - \frac{1}{2\rho}\|z^{k+1}\|^2 \notag \\
	&\mathop{\geq}\limits^{(\ref{z_bound})} \lambda \|u^{k+1}_+\|_0 + \frac{\rho}{2} \| Aw^{k+1} + b - u^{k+1} + \frac{1}{\rho}z^{k+1} \|^2 + \frac{1}{2} f(w^{k+1}) \notag \\
	&~~~ + \frac{1}{2} f(w^{k+1}) - \frac{3}{2\rho\gamma^2} \| \nabla f(w^{k+1}) \|^2 + ( \tau + \frac{\beta}{2} - \frac{3\mu^2}{2\rho\gamma^2} ) \| w^{k+1} - w^k \|^2 + (\eta - \frac{3}{2\rho\gamma^2})\epsilon_{k}^2 \notag \\
	& = \lambda \|u^{k+1}_+\|_0 + \frac{\rho}{2} \| Aw^{k+1} + b - u^{k+1} + \frac{1}{\rho}z^{k+1} \|^2 + \frac{1}{2} f(w^{k+1}) + ( \frac{1}{4l_f} - \frac{3}{2\rho \gamma^2} ) \| \nabla f(w^{k+1}) \|^2  \notag \\
	&~~~ + \frac{1}{2} f(w^{k+1}) - \frac{1}{4l_f} \| \nabla f(w^{k+1}) \|^2 + ( \tau + \frac{\beta}{2} - \frac{3\mu^2}{2\rho\gamma^2} ) \| w^{k+1} - w^k \|^2 + (\eta - \frac{3}{2\rho\gamma^2})\epsilon_{k}^2 \notag \\
	&\mathop{\geq}\limits^{(\ref{f-nabla_f})} \lambda \|u^{k+1}_+\|_0 + \frac{\rho}{2} \| Aw^{k+1} + b - u^{k+1} + \frac{1}{\rho}z^{k+1} \|^2 + \frac{1}{2} f(w^{k+1}) + \frac{1}{2} f( w^{k+1} - \frac{1}{l_f} \nabla f( w^{k+1} ) ) \notag \\
	&~~~+ ( \frac{1}{4l_f} - \frac{3}{2\rho \gamma^2} ) \| \nabla f(w^{k+1}) \|^2 + ( \tau + \frac{\beta}{2} - \frac{\mu^2}{\rho\gamma^2} ) \| w^{k+1} - w^k \|^2 + (\eta - \frac{1}{\rho\gamma^2})\epsilon_{k}^2 \notag \\
	&\geq \lambda \|u^{k+1}_+\|_0 + \frac{\rho}{2} \| Aw^{k+1} + b - u^{k+1} + \frac{1}{\rho}z^{k+1} \|^2 + \frac{1}{2} f(w^{k+1}) + \frac{1}{2}\inf\limits_{w \in \mathbb{R}^p} f(w)  \notag \\
	&~~~  + ( \tau + \frac{\beta}{2} - \frac{\mu^2}{\rho\gamma^2} ) \| w^{k+1} - w^k \|^2 + (\eta - \frac{1}{\rho\gamma^2})\epsilon_{k}^2 + ( \frac{1}{4l_f} - \frac{3}{2\rho \gamma^2} ) \| \nabla f(w^{k+1}) \|^2. \label{bound1}
\end{align}
Then we have
\begin{align}
	&\mathcal{V}_1 + \eta  \epsilon_0^2 - \frac{1}{2}\inf\limits_{w \in \mathbb{R}^p} f(w) \notag \\
	\geq& \lambda \|u^{k+1}_+\|_0 + \frac{\rho}{2} \| Aw^{k+1} + b - u^{k+1} + \frac{1}{\rho}z^{k+1} \|^2  + (\eta - \frac{3}{2\rho\gamma^2})\epsilon_{k}^2 + \frac{1}{2} f(w^{k+1}) \notag \\
	& + ( \tau + \frac{\beta}{2} - \frac{3\mu^2}{2\rho\gamma^2} ) \| w^{k+1} - w^k \|^2 + ( \frac{1}{4l_f} - \frac{3}{2\rho \gamma^2} ) \| \nabla f(w^{k+1}) \|^2. \label{upbound1}
\end{align}
From (\ref{para_mu}), we can justify 
\begin{equation}\notag \label{>0}
	\eta - \frac{3}{2\rho\gamma^2}  > 0, ~~~ \tau + \frac{\beta}{2} - \frac{3\mu^2}{2\rho\gamma^2} > 0, ~~~ \frac{1}{4l_f} - \frac{3}{2\rho \gamma^2} > 0.
\end{equation}
By the coercivity of $f$ and (\ref{upbound1}), we know that $\{ w^{k} \}_{k \in \mathbb{N}}$ is bounded. Since $v^k = w^{k}$, we know that $\{ v_k \}_{k \in \mathbb{N}}$ is bounded.
It follows from (\ref{z_bound1}) and $\epsilon_k < \epsilon_0$ that
\begin{equation}
	\|z^{k+1}\| \leq \frac{1}{\gamma} \| \nabla f(w^{k+1}) \| + \frac{\mu}{\gamma}(\|w^{k+1}\| + \|w^k\|) + \frac{\epsilon_0}{\gamma},
\end{equation}
Since $\{ w^{k} \}_{k \in \mathbb{N}}$ is bounded and $f$ is continuously differentiable, $\{ \| \nabla f(w^{k+1}) \| \}_{k \in \mathbb{N}}$ is bounded and hence $\{z^k\}_{k\in{\mathbb{N}}}$ is bounded. From (\ref{eq2.2}), we have
\begin{equation} \notag
	\|u^{k+1}\| \leq \| A \| \|w^{k+1}\| + \|b\| + \frac{1}{\rho}( \|z^{k+1} \| + \|z^k\| ).
\end{equation}
This implies that $\{u^k\}_{k\in{\mathbb{N}}}$ is bounded. Overall, the sequence $\{ (w^k,u^k,z^k,v^k) \}_{k\in \mathbb{N}}$ generated by IALM is bounded. \hfill\Halmos
\endproof

The following lemma can be immediately deduced from Lemmas \ref{lemma3.1} and \ref{seq_bounded}.

\begin{lemma} \label{cor3.2}
	Under the premises of Lemma \ref{lemma3.1}, we have
	\begin{equation} \label{eq3.7}
		\lim_{k \to \infty} \| w^{k+1} - w^k \| = 0,~ \lim_{k \to \infty} \| u^{k+1} - u^k \| = 0,~
		\lim_{k \to \infty} \| z^{k+1} - z^k \| = 0,~\lim_{k \to \infty} \| v^{k+1} - v^k \| = 0.
	\end{equation}
\end{lemma}
\proof{Proof.} The Lemma \ref{lemma3.1} indicates that $\{ \mathcal{V}_k + \eta \epsilon_{k-1}^2 \}_{k\in \mathbb{N^+}}$ is a decreasing sequence. According to the definition of Lyapunov function, we have
\begin{align}
	\mathcal{V}_k + \eta \epsilon_{k-1}^2 
	=& f(w^k) + \lambda \| u^k_+ \|_0 + \frac{\rho}{2} \| Aw^k + b - u^k +\frac{1}{\rho}z^k \|^2 \label{eq3.9}  + \frac{\beta}{2} \| w^k - v^{k-1} \|^2 + \eta \epsilon_{k-1}^2 - \frac{1}{2\rho}\|z^k\|^2 . \notag
\end{align} 
Notice that $f$ is bounded below and $\{z^k\}_{k\in \mathbb{N}}$ is bounded. We can conclude that $\{ \mathcal{V}_k + \eta \epsilon_{k-1}^2 \}_{k\in \mathbb{N^+}}$ is bounded below. Then $\{ \mathcal{V}_k + \eta \epsilon_{k-1}^2 \}_{k\in \mathbb{N^+}}$ is a decreasing and bounded sequence. Hence, it converges to a finite limit, say $\mathcal{V}_*$,
\begin{equation} \notag \label{eq3.10}
	\lim_{k \to \infty} \mathcal{V}_k + \eta \epsilon_{k-1}^2 = \mathcal{V}_*.
\end{equation}
Taking limit on both sides of (\ref{lyap_des}) implies
\begin{equation} \label{wk+1-k}
	\lim_{k \to \infty} \| w^{k+1} - w^k \| = 0. \notag
\end{equation}
Then $ \lim_{k \to \infty} \| v^{k+1} - v^k \| =  \lim_{k \to \infty} \| w^{k+1} - w^{k} \| = 0.$
Moreover, taking limit as $k \to \infty$ in (\ref{eq3.5}), we have $\lim_{k \to \infty} \| z^{k+1} - z^k \| = 0$.
From (\ref{eq2.2}), we have 
\begin{equation} \notag \label{eq3.12}
	u^{k+1} - u^k = \frac{1}{\rho}( z^k - z^{k-1} ) - \frac{1}{\rho}( z^{k+1} - z^{k} ) + A( w^{k+1} - w^k ),
\end{equation}
which means that
\begin{equation} \label{eq3.13}
	\|u^{k+1} - u^k \| \leq \frac{1}{\rho}\|z^k - z^{k-1}\| + \frac{1}{\rho}\| z^{k+1} - z^{k} \| + \|A\| \| w^{k+1} - w^k \|.
\end{equation}
Then passing to limit as $k \to \infty$ in (\ref{eq3.13}) verifies that $\lim_{k \to \infty} \| u^{k+1} - u^k \| = 0$ and we complete the proof. \hfill\Halmos
\endproof

Now we are ready to show our first convergence result.
\begin{theorem}[Subsequence Convergence] \label{sub_con}Suppose that Assumptions \ref{bbs}, \ref{ass_fullrr} and \ref{assum3.1} hold and parameters are chosen as (\ref{para_mu}) and (\ref{epsilon_val}). If $\{ (w^k;u^k;z^k;v^k) \}_{k\in \mathbb{N}}$ is a sequence generated by IALM satisfying stopping criteria (\ref{stopping1}) and $(w^*;u^*;z^*;v^*)$ is an accumulation point of $\{ (w^k;u^k;z^k;v^k) \}_{k\in \mathbb{N}}$. Then $v^* = w^*$ and $(w^*;u^*;z^*)$ is a {\rm P}-stationary triplet of (\ref{eq2}).
\end{theorem}
\proof{Proof.}
The sequence $\{ (w^k;u^k;z^k) \}_{k\in \mathbb{N}}$ is bounded and therefore there exists a subsequence $\{ (w^{k};u^{k};z^{k}) \}_{k\in \mathcal{K}}$ converges to $(w^*,u^*,z^*)$. Moreover, $\lim_{k \to \infty, k \in \mathcal{K}} v^k = \lim_{k \to \infty, k \in \mathcal{K}} w^k = w^* $. From (\ref{crit1}), (\ref{nabla_h}) and (\ref{eq2.2}), we can obtain
\begin{align}
	{\rm dist}( u^{k}, {\rm Prox_{\alpha\lambda \|(\cdot)_+\|_0}} ( u^{k} + \alpha z^k ) ) \leq \epsilon_k.
\end{align}
According to the proximal behavior in \citet[Theorem 1.25]{RockWets98}, we know that ${\rm Prox_{\alpha\lambda \|(\cdot)_+\|_0}} ( u^{k} + \alpha z^{k} )$ is a nonempty compact set. Then taking (\ref{crit1}) into account, there exists a point $\overline{u}^{k} \in {\rm Prox_{\alpha\lambda \|(\cdot)_+\|_0}} ( u^{k} + \alpha z^{k} )$ such that 
\begin{equation} \notag \label{ud_dis}
	\| u^{k} - \overline{u}^{k} \| \leq \epsilon_{k-1},~~\forall k = 1,2,\cdots.
\end{equation}
Since $\lim_{k \to \infty} \epsilon_k = 0$, we have $ \lim_{k \to \infty} u^{k} - \overline{u}^{k} = 0$, which implies that
\begin{equation} \notag
	\lim_{k\in \mathcal{K},k \to \infty}  \overline{u}^{k} = \lim_{k\in \mathcal{K},k \to \infty} u^{k} - (u^{k} - \overline{u}^{k}) = u^*. 
\end{equation}
It follows from the proximal behavior (see \citet[Theorem 1.25]{RockWets98}) that 
\begin{equation} \label{sub_Psta}
	u^* \in {\rm Prox_{\alpha\lambda \|(\cdot)_+\|_0}} ( u^{*} + \alpha z^{*} ).
\end{equation}
From (\ref{crit1}) and (\ref{eq2.2}), we have
\begin{align}
	& \| \nabla f(w^k) + \mu ( w^k - v^{k-1} ) + A^\top z^{k} \| = \| \nabla_w h_{k-1}(w^k,u^k) \| \leq \epsilon_{k-1}, \label{subw} \\
	& z^k - z^{k-1} = \rho ( Aw^k + b - u^k ). \label{subz}
\end{align}	
Notice that $v^{k-1} = w^{k-1}$ and (\ref{eq3.7}) holds. Taking $k \in \mathcal{K}$, $k \to \infty$ on both side of (\ref{subw}) and (\ref{subz}), we have	
\begin{equation} \notag \label{sub_Psta1}
	\nabla f(w^*) + A^\top z^* = 0, ~~
	Aw^* + b = u^* .
\end{equation}
Combining with (\ref{sub_Psta}), we can conclude that $ (w^*,u^*,z^*) $ is a P-stationary triplet of (\ref{eq2}).\hfill\Halmos
\endproof

\begin{remark}
	It is noteworthy that the optimality conditions of (\ref{eq2}) and (\ref{exact_penalty}) are established by proximal-type stationary points and the stopping criteria (\ref{crit1}) is defined based on a proximal operator.  From the proof of Theorem \ref{pam_converge} (iv) and Theorem \ref{sub_con}, we can see that the proximal behavior is essential to the subsequence convergence.
\end{remark}

We establish the whole sequence convergence under the following assumption, which is stronger than Assumption \ref{assum3.1}.

\begin{assumption} \label{ass_str_con}
	The function $f$ is strongly convex.
\end{assumption}


\begin{theorem}[Whole Sequence Convergence]\label{whole_con}Suppose that Assumptions \ref{bbs}, \ref{ass_fullrr} and \ref{ass_str_con} hold, and parameters are chosen as (\ref{para_mu}) and (\ref{epsilon_val}). If $\{ (w^k;u^k;z^k;v^k) \}_{k\in \mathbb{N}}$ is a sequence generated by IALM satisfying stopping criteria (\ref{stopping1}), then $\lim_{k \to \infty} (w^k;u^k;z^k;v^k) = (w^*;u^*;z^*;w^*)$,  $(w^*;u^*;z^*)$ is a {\rm P}-stationary triplet of (\ref{eq2}) and $(w^*;u^*)$ is a local minimizer of (\ref{eq2}).
\end{theorem}
\proof{Proof.}
Let us first prove that the number of {\rm P}-stationary triplet of (\ref{eq2}) is finite. Suppose that $(w^*,u^*,z^*)$ is a {\rm P}-stationary triplet of (\ref{eq2}), then it satisfies (\ref{P-stat}). Since the given $u^*$ and $z^*$ satisfy the second formula in (\ref{P-stat}), similar to the analysis in Property \ref{rela_cr_p}, we can deduce that $u^*$ and $z^*$ satisfy
\begin{equation}  \label{prox_uz_val}
	\left\{ \begin{aligned}
		& z^*_i = 0  , &&u^*_i \in (-\infty,0] \cup [\sqrt{2\lambda\alpha}, \infty), \\
		& u^*_i = 0, &&  z^*_i \in [0,\sqrt{2\lambda/\alpha}],
	\end{aligned}  \right. ~\forall i \in [n].
\end{equation}
Denote $ S^*_- := \{ i \in [n]: ( Aw^* + b )_i \leq 0 \} $. Let us consider the following convex programming:
\begin{equation} \label{convex_pro1}
	\begin{aligned}
		&\min_{w \in \mathbb{R}^p} &&f(w)  \\
		&~~{\rm s.t.} && (Aw + b )_{S^*_-} \leq 0.
	\end{aligned}
\end{equation}
Note that the objective function $f
(w)$ is strongly convex and the constrains are linear inequalities. Therefore, (\ref{convex_pro1}) has unique global minimizer. Moreover, the KKT condition of (\ref{convex_pro1}) can be written as 
\begin{equation} \label{kkt_sys1}
	\left\{  \begin{aligned}
		& \nabla f(w) + A_{S^*_-}^\top z_{S^*_-} = 0, \\
		& z_{S^*_-} \geq 0,~ u_{S^*_-} \leq 0,~ \langle z_{S^*_-},~ u_{S^*_-}  \rangle = 0, \\
		& (Aw + b - u)_{S^*_-} = 0,
	\end{aligned}  \right.
\end{equation}
where $z_{S^*_-} \in \mathbb{R}^{|S^*_-|}$ is the Lagrange multiplier. It follows from (\ref{P-stat}) and (\ref{prox_uz_val}) that $(w^*,u^*,z^*)$ satisfies (\ref{kkt_sys1}) and hence $w^*$ is the unique global minimizer of (\ref{convex_pro1}). Moreover, $	u^* = Aw^* + b$ and $z^* = -(AA^\top)^{-1}\nabla f(w^*)$ yields that $u^*$ and $z^*$ can be identified if $w^*$ is fixed. Notice that there are only a finite number of $S^*_-$ (no more than $2^n$). Therefore, the number of P-stationary triplet for (\ref{eq2}) is finite.

From the conclusion of Theorem \ref{sub_con}, we know that each accumulation point of sequence $\{ (w^k;u^k;z^k) \}_{k\in \mathbb{N}}$ is a P-stationary triplet of (\ref{eq2}). Hence, there is a finite number of accumulation points for $\{ (w^k;u^k;z^k) \}_{k\in \mathbb{N}}$, which means that each accumulation point of $\{ (w^k;u^k;z^k) \}_{k\in \mathbb{N}}$ is isolated. Considering that (\ref{eq3.7}) holds true, it follows from \citet{kanzow1999qp} that $\{ (w^k;u^k;z^k) \}_{k\in \mathbb{N}}$ must converges to a P-stationary triplet $(w^*;u^*;z^*)$. Moreover, $\lim_{k \to \infty} v^k = \lim_{k \to \infty} w^{k} = w^* $. Finally, Lemma \ref{opt_orig} (ii) implies that $(w^*;u^*)$ is a local minimizer of (\ref{eq2}). \hfill \Halmos
\endproof


\section{Applications.} \label{application}
In this section, we will show some examples of 0/1-COP which can be solved by our IALM with global convergence.
\subsection{Support Vector Machine.} The classical support vector machine (SVM) proposed by \citet{cortes1995support} is a useful tool for binary classification. It aims to construct two parallel support hyperplanes separating two classes meanwhile the width between the two parallel hyperplanes is maximized. Given a training set $\{ (x_i,y_i) : i \in [n] \}$ with the $i$th sample $x_i \in \mathbb{R}^{p}$, where $x_{ip} = 1$, and the $i$th class label $y_i \in \{ -1,1 \}$, we denote $y :=(y_1,\cdots,y_n)^\top \in \mathbb{R}^n$ and $X:=(x_1,x_2,\cdots,x_n)^\top \in \mathbb{R}^{n \times p}$. Referring to the zero-one loss term in \citet{wang2021support} and the structural risk term in \citet{lee2001ssvm,mangasarian1999successive}, we proposed the following 0/1-COP for SVM:
\begin{align} \label{0/1-SVM}
	\min_{w \in \mathbb{R}^p} \frac{1}{2} \| w \|^2 + \lambda \|[-(( y\cdot \textbf{1}^\top_p) \odot X)w + \textbf{1}_n ]_+\|_0.
\end{align}
The above problem can be recovered from (\ref{eq1}) by taking
\begin{equation} \label{SVM_set}
	f(w) :=  \frac{1}{2} \| w \|^2,~ A := -( y\cdot \textbf{1}^\top_p) \odot X,~b := \textbf{1}_n,~\lambda > 0.
\end{equation}
In the setting of (\ref{SVM_set}), $f$ is a strongly convex function, if we further assume that $A$ is full row rank and choose parameters satisfying (\ref{para_mu}) and (\ref{epsilon_val}), then the whole sequence generated by our IALM converges to a local minimizer of (\ref{0/1-SVM}).

\subsection{Twin Support Vector Machine.}
Differing with the classical SVM, the twin support vector machine (TSVM) proposed by \citet{khemchandani2007twin}  seeks two nonparallel hyperplanes such that each one is closest to one class of data samples meanwhile far away from the other class. It is implemented by solving two smaller composite optimization problem instead of a single large one, which promotes learning speed compared with classical SVM. 

Given the positive training instances $\{ x^{(1)}_i \in \mathbb{R}^{p}: x^{(1)}_{i,p} = 1, i \in [n_1] \}$ and negative training instances $\{ x^{(-1)}_i \in \mathbb{R}^{p} : x^{(-1)}_{ip} = 1, i \in [n_{-1}] \}$, we denote 
matrices $X^{(1)} := ({x}^{(1)}_1,\cdots,x^{(1)}_{n_1})^\top\in \mathbb{R}^{n_1 \times p}$ and $ X^{(-1)} := (x^{(-1)}_1,\cdots,x^{(-1)}_{n_{-1}})^\top \in \mathbb{R}^{n_{-1} \times p} $. In particular, the TSVM with $\ell_2$ regularization (see e.g. \citet{shao2011improvements}) can be represented as the following 0/1-COP:
\begin{align}
	&\min_{w^{(1)} \in \mathbb{R}^p} \frac{1}{2} \| w^{(1)} \|^2 + \frac{\lambda_1}{2} \| X^{(1)} w^{(1)} \|^2 + \lambda_2 \| (\textbf{1}_{n_{-1}} + X^{(-1)} w^{(1)})_+ \|_0, \label{tsvm1}~\lambda_1,\lambda_2 > 0,  \\
	&\min_{w^{(-1)} \in \mathbb{R}^p} \frac{1}{2} \| w^{(-1)} \|^2 + \frac{\lambda_3}{2} \| X^{(-1)} w^{(-1)} \|^2 + \lambda_4 \| (\textbf{1}_{n_{1}} - X^{(1)} w^{(-1)})_+ \|_0,~\lambda_3,\lambda_4 > 0. \label{tsvm-1}
\end{align}
Obviously, (\ref{tsvm1}) is covered by (\ref{eq1}) with
\begin{equation} \notag
	f(w) := \frac{1}{2} \| w \|^2 + \frac{\lambda_1}{2} \| X^{(1)} w \|^2,~ A := X^{(-1)},~ b := \textbf{1}_{n_{-1}},~\lambda:= \lambda_2 > 0.
\end{equation}
Similarly, (\ref{tsvm-1}) can be recovered from (\ref{eq1}) with
\begin{equation} \notag
	f(w) := \frac{1}{2} \| w \|^2 + \frac{\lambda_3}{2} \| X^{(-1)} w \|^2,~ A := - X^{(1)},~ b := \textbf{1}_{n_{1}},~\lambda:= \lambda_4 > 0.
\end{equation}
For each of (\ref{tsvm1}) and (\ref{tsvm-1}), $f$ is a strongly convex function. If we assume that $A$ is full row rank and the parameters are set as (\ref{para_mu}) and (\ref{epsilon_val}), then the whole sequence generated by our IALM converges to a local minimizer of the corresponding problems.

\subsection{Multi-label Classification.} In multi-label classification (MLC), each instance is associated with a subset of class labels instead of a single class label as in traditional classification (see e.g. \citet{zhang2013review}). Particularly, the hamming loss is an important metric for evaluating a multi-label classifier. Now let us consider a MLC problem with $m$ classes. Given the training instances $x_i \in \mathbb{R}^p$ with $x_{i,p} = 1$ for all $ i \in [n]$ and the class labels $y_i  \in \{ 1,-1 \}^m$ for all $i \in [n]$ with the relevant (resp. irrelevant) classes being $1$ (resp. $-1$), we denote $X := ( x_1,\cdots,x_n )^\top \in \mathbb{R}^{n \times p}$ and $Y := (y_1, \cdots, y_n)^\top \in \mathbb{R}^{n \times m} $. For a multi-label classifier $\mathcal{C}(\cdot): \mathbb{R}^p \to \mathbb{R}^m$, the empirical hamming loss could be represented as
\begin{align} \notag
	\ell_H (\mathcal{C},X,Y) := \frac{1}{2nm} \sum_{i = 1}^n \| \mathcal{C}(x_i) - y_i \|_1.
\end{align}  
The linear binary relevance method (see \citet{boutell2004learning}), which has been attracted considerable attention in MLC field, assumes that the multi-label classifier has the intuitive form $\mathcal{C}(x) = ( {\rm sgn}({w^{(1)}}^\top x), \cdots, {\rm sgn}({w^{(m)}}^\top x) )$, where $w^{(k)} \in \mathbb{R}^p$, $\forall k \in [m]$ and the function ${\rm sgn}: \mathbb{R} \to \mathbb{R}$ is defined as 
\begin{align} \notag
	{\rm sgn}(t) := \left\{ \begin{aligned}
		& ~~~1, &&t > 0, \\
		&-1, &&t \leq 0.
	\end{aligned} \right.
\end{align}
Under this assumption, minimizing the hamming loss could be decomposed into the following optimization problems:
\begin{align} \label{br_mlc} 
	\min_{w^{(k)} \in \mathbb{R}^p} \frac{1}{2} \| w^{(k)} \|^2 + \lambda \|[-(( y^{(k)}\cdot \textbf{1}^\top_p) \odot X)w^{(k)} + \textbf{1}_n ]_+\|_0, ~ k \in [m],
\end{align}
where $y^{(k)} \in \mathbb{R}^n$ is the $k$th column of $Y$. The $k$th optimization problem of (\ref{br_mlc}) could be seen as a special case of (\ref{eq1}) with
\begin{align} \notag
	f(w) :=  \frac{1}{2} \| w \|^2,~ A := -( y^{(k)}\cdot \textbf{1}^\top_p) \odot X,~b := \textbf{1}_n,~\lambda > 0.
\end{align}
Since $f$ is strongly convex, if $A$ is full row rank and parameters are set as (\ref{para_mu}) and (\ref{epsilon_val}), then the whole sequence generated by our IALM converges to a local minimizer of the corresponding problems.

\subsection{Ridge Regression with Maximum Rank Correlation.} 
\citet{han1987non} proposed the maximum rank correlation (MRC) to estimate the coefficients of generalized regression model. The MRC estimator is asymptotically normal (see e.g. \citet{fan2020rank}) and robust to heavy tailed error distribution (see e.g. \citet{li2021central}). 
Given observations $x_i \in \mathbb{R}^p$ and $y_i \in \mathbb{R}$ for $i \in [n]$, the MRC method aims to maximize the following measure
\begin{align} \notag
	\mathcal{S}(w):= \frac{1}{n(n-1)} \sum_{i \neq j} \mathbb{I} (y_i < y_j) \mathbb{I}( w^\top x_i < w^\top x_j ),
\end{align} 
where $\mathbb{I}(\cdot)$ is an indicator function,
\begin{align}
	\mathbb{I} ( \cdot ) := \left\{ \begin{aligned}
		& 1, &&{\rm if~(\cdot)~is~true}, \\
		& 0, &&{\rm otherwise}.
	\end{aligned} \right.
\end{align}
It is not hard to see that $\mathcal{S}(w) \in [0,1]$. Without loss of generality, suppose that the $n$ observations of response variables could be ranked as $y_1 \leq y_3 \leq \cdots \leq y_n$. Then $w^\top x_1 < w^\top x_2 < \cdots < w^\top x_n $ ensures $\mathcal{S}(w) = 1$. Therefore, to achieve MRC estimation, we only need to consider minimizing the following term 
\begin{align}
	\mathcal{R}(w) := \| (\textbf{1}_{n-1} \cdot \xi + BXw)_+ \|_0,
\end{align}
where 
\begin{align} \notag
	X := (x_1,\cdots,x_n)^\top \in \mathbb{R}^{n \times p},~	B := \left[ \begin{array}{ccccc}
		1 &-1 & & &\\
		& 1 &-1 & & \\
		&   & \ddots &\ddots & \\
		& & &1 &-1
	\end{array}  \right] \in \mathbb{R}^{(n-1)\times n},~\xi > 0.
\end{align}
Based on this idea, to make a trade-off between the data fitting and MRC, we give the following optimization for the ridge regression with MRC:
\begin{align} \label{rr_mrc}
	\min_{w \in \mathbb{R}^p}	\frac{1}{2} \| Xw - y \|^2 + \frac{\lambda_1}{2} \| w \|^2 + \lambda_2 \mathcal{R}(w),~ \lambda_1, \lambda_2 > 0,
\end{align}
where $y := ( y_1,\cdots, y_n )^\top \in \mathbb{R}^n$. Obviously, it is a special case of (\ref{eq1}) with
\begin{align} \notag
	f(w):= \frac{1}{2} \| Xw - y \|^2 + \frac{\lambda_1}{2} \| w \|^2,~ A = BX,~b:=\textbf{1}_{n-1}\cdot \xi,~\lambda : = \lambda_2 > 0.
\end{align}
Since $f$ is strongly convex, if $A$ is full row rank and parameters are set as (\ref{para_mu}) and (\ref{epsilon_val}), then the whole sequence generated by our IALM converges to a local minimizer of (\ref{rr_mrc}).

\section{Conclusion.}
In this paper, we studied the property of Lyapunov function and proposed an inexact augmented Lagrangian method (IALM) with global convergence for solving 0/1-COP. We first utilized a proximal-type stationarity to derive the minimum and the (strongly) exact penalty property of the Lyapunov function with a fixed Lagrange multiplier. Then we design an IALM for solving 0/1-COP and the subproblems generated by IALM are solved by a zero-one Bregman alternating linearized minimization (0/1-BALM) method. We proved the global convergence of our IALM equipped with 0/1-BALM under some suitable assumptions. Finally, we applied our IALM for solving SVM, MLC and MRC, and showed global convergence of IALM with the full row rank of data matrix in these cases. A issue that needs further consideration is how to derive the local convergence rate of our IALM. We leave this issue to be explored in future research.



%
%
%

\section*{Acknowledgments.}
This work was supported by the National Natural Science Foundation of China (12131004, 11971052) and Beijing Natural Science Foundation (Z190002).


\bibliographystyle{informs2014} 
\bibliography{references} 


\end{document}